\documentclass[a4paper,11pt,pdftex,bibtotocnumbered, headsepline,normalheadings%, draft%twoside
]%{scrreprt}  %draft in der eckigen Klammer zeigt ueberstehende Zeilen 
{article}
\usepackage{subfigure,fancyhdr,pstricks}
\usepackage{amssymb,amsmath,amsfonts,amsthm,stmaryrd}
\usepackage[T1]{fontenc}
\usepackage{tensor}
\usepackage{epsfig}
\usepackage{graphicx}
\usepackage{makeidx}
\usepackage[enableskew]{youngtab}
\usepackage{rotating}
\usepackage[all, knot]{xy}
\xyoption{poly}
\xyoption{arc}
\usepackage{lscape}
\usepackage[colorlinks=true,linkcolor=black,pdftex]{hyperref}
\usepackage[hyperpageref]{backref} 
\usepackage{faktor}
\usepackage[nice]{nicefrac}
\usepackage{calc}
\usepackage{niceframe}
\usepackage{tikz}
\usetikzlibrary{shapes,arrows}

\tikzstyle{decision} = [diamond, draw, text width=4.5em, text badly centered, node distance=3cm, inner sep=0pt]
\tikzstyle{block} = [rectangle, draw, text width=5em, text centered, rounded corners, minimum height=4em]

\theoremstyle{plain}

\newtheorem{satz}{Theorem}[section]

\newtheorem{coro}{Corollary}[section]
\newtheorem{lemma}{Lemma}

\theoremstyle{definition}
\newtheorem{defi}{Definition}
\theoremstyle{definition}

\newtheorem*{exa}{Example}
\newtheorem*{rem}{Remark}

\theoremstyle{remark}

\DeclareMathOperator{\id}{id}

\DeclareMathOperator{\Res}{Res}

\DeclareMathOperator{\Ext}{\rm Ext}
\DeclareMathOperator{\Simp}{\rm Simp}

\DeclareMathOperator{\Hom}{\rm Hom}
\DeclareMathOperator{\End}{\rm End}

\DeclareMathOperator{\inj}{\hookrightarrow}

\newcommand{\oTo}{\xymatrix{ \ar@{^{(}->}[r]|{\mathbf{O}}& }}
\newcommand{\cTo}{\xymatrix{ \ar@{^{(}->}[r]|{\mathbf{|}}& }}
\newcommand{\coTo}{\xymatrix{ \ar@{^{(}->}[r]|{\mathbf{O}}|{\mathbf{|}}& }}
\DeclareMathOperator{\surj}{\twoheadrightarrow }

\newcommand{\tits}[2]{\langle #1, #2 \rangle}

\DeclareMathOperator{\Ker}{ker}

\DeclareMathOperator{\Ind}{Ind}
\DeclareMathOperator{\Stab}{Stab}
\DeclareMathOperator{\lin}{lin}
\DeclareMathOperator{\Mod}{mod}
\DeclareMathOperator{\Lie}{Lie}
\DeclareMathOperator{\codim}{codim}
\DeclareMathOperator{\res}{Res}

\newcommand{\Gl}{\mathbf{Gl}}

\newcommand{\Sl}{\mathbf{Sl}}

\newcommand{\si}{\sigma}

\newcommand{\la}{\lambda}

\newcommand{\al}{\alpha}

\newcommand{\D}{\mathbb{D}}

\newcommand{\HH}{\mathbb{H}}

\newcommand{\N}{\mathbb{N}}

\newcommand{\Z}{\mathbb{Z}}
\newcommand{\Q}{\mathbb{Q}}

\renewcommand{\S}{\mathbb{S}}
\newcommand{\C}{\mathbb{C}}

\newcommand{\mcA}{\mathcal{A}}
\newcommand{\mcB}{\mathcal{B}}

\newcommand{\mcH}{\mathcal{H}}

\newcommand{\mcL}{\mathcal{L}}

\newcommand{\mcN}{\mathcal{N}}
\newcommand{\mcO}{\mathcal{O}}
\newcommand{\mcP}{\mathcal{P}}

\newcommand{\mcS}{\mathcal{S}}

\newcommand{\mcU}{\mathcal{U}}

\newcommand{\Rd}{\rm R_Q(\underline{d})}

\newcommand{\Fl}{\rm Fl}

\newcommand{\dd}{\underline{d}}

\newcommand{\ee}{\underline{e}}
\newcommand{\ff}{\underline{f}}
\newcommand{\ddd}{\underline{\bold d}}

\newcommand{\Gd}{\mathbf{Gl}_{\underline{d }}}

\begin{document}

%\include{0-chapter-Notation}

%-----------------------------------------------------------------------------------------------------
\newpage
%-----------------------------------------------------------------------------------------------------
%%%%%%%%%%%%%%%%%%%%%%%%--0--%%%%%%%%%%%%%%%%%%%%%%%%%%%%%%%%%%%%%%%%%%%%%%%%%%%
\section*{A Survey on Springer theory (over $\C$)}
\begin{center} \today \end{center}

\begin{abstract} A \textbf{Springer map} is for us a union of collapsings of (complex) 
homogeneous vector bundles and a \textbf{Steinberg variety} is just the cartesian product of a Springer map with itself. 
Ginzberg constructed on the (equivariant) Borel-Moore homology and on the (equivariant) $K$-theory of a Steinberg variety a convolution product making it an associative algebra, we call this a \textbf{Steinberg algebra}. The decomposition theorem for perverse sheaves gives the indecomposable, projective graded 
modules over the Steinberg algebra. Also Ginzberg's convolution yields a module structure on the respective homology groups of the fibres under the Springer maps, which we call \textbf{Springer fibre modules}. In short, 
\emph{for us} a \textbf{Springer theory} is the study of a Steinberg algebra together with its graded modules. \\
We give two examples: Classical Springer theory and quiver-graded Springer theory.\\
%Then we turn back to the general setup and explore some basic properties. 
%Zwei neue Themen miteinbauen - Carrell,deConcini-Procesi-Thm and cellularity, cellulary stratified and quasi-hereditary?
 
\begin{itemize}
\item[(1)] Definitions and basic properties.
\item[(2)] Examples
\begin{itemize}
\item[(a)] Classical Springer Theory.
%\item[(b)] Exotic Springer Theory.
\item[(b)] Quiver-graded Springer Theory.
%Here we also study induced modules.
\end{itemize}
\item[(3)] We discuss literature on the two examples.   
\end{itemize}
\end{abstract}

\subsection*{Definition of a Springer Theory}

Roughly, following the introduction of Chriss and Ginzburg's book (\cite{CG})\footnote{We take a more general approach, what usually is considered as Springer theory you find in the example \emph{classical Springer theory}. Nevertheless, our approach is still only a special case of \cite{CG}, chapter 8.}, 
Springer Theory is a uniform geometric construction for a wide class of (non-commutative) algebras together with families of modules 
over these algebras. Examples include 
\begin{itemize}
\item[(1)] Group algebras of Weyl groups together with their irreducible representations,
\item[(2)] affine Hecke algebras together with their standard modules and irreducible representations,
\item[(3)] Hecke algebras with unequal parameters, 
\item[(4)] KLR-algebras ($=$ Quiver Hecke algebras) 
\item[(5)] Quiver Schur algebras
\end{itemize}
To understand this construction, recall  
for any algebraic group $G$ and closed subgroup $P$ (over $\C$) we call the principal bundles $G\to G/P$ \textbf{homogeneous}. 
For any $P$-variety $F$ given we have the \textbf{associated bundle} defined by the quotient 
\[G\times^P F:=G\times F/\sim \; , \quad (g,f)\sim (g^{\prime},f^{\prime})\colon \;\iff \text{ there is } p\in P \colon (g,f)=(g^{\prime}p,p^{-1}f^{\prime})\]
 and $G\times^P F\to G/P, (g,f)\mapsto gP$.  
Given a representation $\rho\colon P\to \Gl(F)$, i.e. a morphism of algebraic groups, we call associated bundles of the form   
$G\times^P F \to G/P$ \textbf{homogeneous vector bundles}. 

\begin{defi}

The uniform geometric construction in all cases is given by the following: 
Given $(G,P_i,V,F_i)_{i\in I}$ with $I$ some finite set,  
\[
\left[
\begin{aligned}
(*) & \; G\text{ a connetcted reductive group with parabolic subgroups }P_i. \\
&\text{We also assume there exists a maximal torus }T\subset G\text{ which is contained in every }P_i.\\ 
(*) & \; V\text{ a finite dimensional $G$-representation, }F_i\subset V\text{ a }P_i\text{-subrepresentation of }V, \; i\in I.%\text{, define }F:=\bigsqcup_{i\in I}F_i. 
\end{aligned}
\right.
\]
We identify $V, F_i$ with the affine spaces having the vector spaces as $\C$-valued points and consider the following morphisms 
of algebraic varieties\footnote{algebraic variety = separated integral scheme of finite type over a field}, let $E_i:= G\times^{P_i}F_i, i\in I$ 

\[
\xymatrix{
 &E:= \bigsqcup_{i\in I}E_i \ar[dl]_{\pi}\ar[dr]^\mu & && &[(g,f_i)] \ar@{|-_{>}}[dl]\ar@{|-_{>}}[dr] & \\
V & & \bigsqcup_{i\in I} G/P_i && gf_i & & gP_i
}
\]
Then, $E\to V\times \bigsqcup_{i\in I} G/P_i, [(g,f_i)]\mapsto (gf_i, gP_i)$ is a closed embedding (see \cite{Sl}, p.25,26), it follows that $\pi$ 
is projective. We call the algebraic correspondence\footnote{two scheme morpisms $\xymatrix{X&Z \ar[l]_p \ar[r]^q &Y}$ are called algebraic correspondence, if $p$ is proper and $q$ is flat} $(E,\pi ,\mu )$  \textbf{Springer triple}\index{Springer triple}, 
the map $\pi$ \textbf{Springer map}\index{Springer map}, its fibres  \textbf{Springer fibres}\index{Springer fibres}. 
Via restriction of $E\to V\times \bigsqcup_{i\in I} G/P_i$ to $\pi^{-1}(x)\to \{x\}\times\bigsqcup_{i\in I}  G/P_i$ one sees that all Springer fibres are via $\mu$ closed subschemes of $\bigsqcup_{i\in I} G/P$.\\

We also have another induced roof-diagramm  
\[
\xymatrix{ & Z:=E\times_V E \ar[dl]_{p}\ar[dr]^m& \\ 
V & & (\bigsqcup_{i\in I}G/P_i) \times(\bigsqcup_{i\in I}G/P_i)  
}
\]
with $p\colon  E\times_V E^\prime \xrightarrow{pr_E}E\xrightarrow{\pi}V$ projective and 
$m\colon  E\times_V E^\prime \xrightarrow{(pr_E,pr_E)}E\times E \xrightarrow{\mu\times \mu}(\bigsqcup_{i\in I}G/P_i) \times(\bigsqcup_{i\in I}G/P_i)$. Observe, by definition 
\[ Z = \bigsqcup_{i,j \in I} Z_{i,j}, \quad Z_{i,j}= E_i \times_V E_j .\] 
We call the roof-diagram $(Z,p,m)$ \textbf{Steinberg triple}\index{Steinberg triple},
the scheme $Z$ \textbf{Steinberg variety}\index{Steinberg variety} (even though as a scheme $Z$ might be neither reduced nor irreducibel). But in 
view of our (co-)homology choice below we only study the underlying reduced scheme and look at its $\C$-valued points endowed with the analytic topology.
\end{defi}

If all parabolic groups $P_i$ are Borel groups, the Steinberg variety $Z$ is an iterated cellular fibration over $\bigsqcup_{i\in I}G/P_i$ (for an apropiate definition of iterated cellular fibration), for the precise statement see the next lemma.  
We choose a (co-)homology theory which can be calculated for spaces with cellular fibration property and which has a localization to the $T$-fixpoint theory.  
Let $H_*^A,\; A\in \{pt, T,G\}$ be \textbf{($A$-equivariant) Borel-Moore homology}. 
We could also choose (equivariant) $K$-theory, but we just give some known results about it.\\ 
There is a natural product $*$ on $H_*^A(Z)$ called convolution product constructed by Chriss and Ginzburg in \cite{CG}. 
\[
\begin{aligned}
*\colon H_*^A(Z) \times H_*^A(Z) & \to H_*^A(Z) \\
   (c_{1,2}, c_{2,3}) &\mapsto c_{1,2}*c_{2,3}:=(q_{1,3})_*(p_{1,2}^*(c_{1,2}) \cap p_{2,3}^*(c_{2,3})) 
\end{aligned}
\] 
where $\cap \colon H_p^A(X) \times H_q^A(Y) \to H_{p+q-2d}^A(X\cap Y)$ is the intersection pairing which is induced by the $\cup$-product in relative singular 
cohomology for $X,Y\subset M$ two $A$-equivariant closed subsets of a $d$-dimensional complex manifold $M$ (cp. \cite{CG}, p.98, (2.6.16)) and
where $p_{a,b} \colon E\times E\times E \to E \times E$ is the projection on the $a,b$-th factors, $q_{a,b}$ is the restriction of $p_{a,b}$ to $E\times_VE\times_VE$. 
It holds 
\[H_p^A(Z_{i,j})* H_q^A(Z_{k,\ell}) \subset \delta_{j,k} H_{p+q-e_k}^A(Z_{i, \ell}), \quad e_k= \dim_{\C}E_k .\]
We call $(H_*^A(Z),*)$ the ($A$-equivariant) \textbf{Steinberg algebra} for $(G,P_i,V,F_i)_{i\in I}$.
%It holds that $H_*^A(Z)$ is a $H_*^A(pt)$-algebra by definition.  

There is a the following identification.  
\begin{satz} \label{Extalgebra}(\cite{CG}, chapter 8) Let $A\in \{pt, T,G\}$ we write  $e_i =\dim_{\C} E_i$. There is an isomorphism of $\C$-algebras 
\[ H_*^A(Z) \to \Ext^*_{D^b_A(V)}(\bigoplus_{i\in I}{(\pi_i)}_*\underline{\C}[e_i], \bigoplus_{i\in I}{(\pi_i)}_*\underline{\C}[e_i]) .\]
If we set 
\[ H_{[p]}^A(Z) := \bigoplus_{i,j \in I} H_{e_i+e_j-p}^A(Z_{i,j})
\]
then $H_{[*]}^A(Z)$ is a graded $H_A^*(pt)$-algebra and the isomorphism is an isomorphism of graded algebras. Furthermore, 
the Verdier duality on $D_A^b(V)$ induces an anti-involution on this algebra.  
\end{satz}
the proof is only given for $A=pt$, but as Varagnolo and Vasserot in \cite{VV} observed, the same proof can be rewritten for the $A$-equivariant case.

\subsubsection*{Convolution Modules (see \cite{CG}, section 2.7)}

Given two subsets $S_{1,2}\subset M_1\times M_2, \;S_{2,3}\subset M_2\times M_3$ the set-theoretic convolution is defined as 
\[ S_{1,2}\circ S_{2,3}:= \{ (m_1,m_3)\mid \exists \;m_2\in M_2 \colon (m_1,m_2)\in S_{1,2}, (m_2,m_3)\in S_{2,3}\} \subset M_1\times M_3 .\]
Now, let $S_{i,j}\subset M_i\times M_j$ be $A$-equivariant locally closed subsets of smooth complex $A$-varieties, let $p_{i,j}\colon M_1\times M_2\times M_3\to M_i\times M_j$ be projection on the $(i,j)$-th factors and assume $q_{1,3}:=p_{1,3}|_{p_{12}^{-1}(S_{1,2})\cap p_{2,3}^{-1}(S_{2,3})} $ is proper. Then we get a map 
\[
\begin{aligned}
*\colon H_p^A(S_{1,2}) \times H_q^A(S_{2,3}) &\to H^A_{p+q-2\dim_{\C} M_2} (S_{1,2}\circ S_{2,3})\\
c_{1,2}*c_{2,3} &:= (q_{1,3})_* (p_{1,2}^*c_{1,2}\cap p_{2,3}^*c_{2,3}).
\end{aligned}
\] 
This way we defined the algebra structure on the Steinberg algebra, 
but it also gives a left module stucture on $H_*^A(S)$ for any $A$-variety $S$ with $Z\circ S=S$ and a right module structure when $ S\circ Z= S$. 
\begin{itemize}
\item[(a)] $M_1=M_2=M_3=E$, embedd $Z=E\times_VE\subset E\times E, E=E\times pt \subset E\times E$, then it holds $Z\circ E =E$. 
If we regrade the Borel-Moore homology (and the Poincare dual $A$-equivariant cohomology) of $E$ as follows 
\[ H_{[p]}^A(E):= \bigoplus_{i\in I} H_{e_i-p}^A(E_i) \quad (=\bigoplus_{i\in I} H_A^{e_i+p}(E_i) =: H_A^{[p]}(E)) \]
then $H_{[*]}^A(E)$ and $H^{[*]}(E)$ carry the structure of a graded left $H_{[*]}^A(Z)$-module. 
\item[(b)] $M_1=M_2=M_3=E$, embedd $E\subset E\times E$ diagonally, then $E\circ E= E$, it holds $H_{(*)}^A(E) =H_A^*(E)$ as graded algebras 
where $H_{(p)}^A(E):= \bigoplus_{i} H_{2e_i-p}^A (E_i)$ and the ring structure on the cohomology is given by the cup product. 
If we take now $Z=E\times_VE\subset E\times E$ then $E\circ Z =Z$ and we get a structure as 
graded left $H_A^*(E)$-module on $H_{[*]}^A(Z)$. 
\item[(c)] $M_1=M_2=M_3=E$, $A=pt$ embedd $Z=E\times_VE\subset E\times E, \pi^{-1}(x)=\pi^{-1}(x)\times pt \subset E\times E$, then it holds $Z\circ \pi^{-1}(x) =E$. 
If we regrade the Borel-Moore homology and singular cohomology of $\pi^{-1}(x)$ as follows 
\[    H_{[p]}(\pi^{-1}(x)):= \bigoplus_{i\in I} H_{e_i-p}(\pi_i^{-1}(x)) \quad  H^{[p]}(\pi^{-1}(x)):= \bigoplus_{i\in I} H^{e_i+p}(\pi_i^{-1}(x)) \]
then $H_{[*]}(\pi^{-1}(x)$ and $H^{[*]}(\pi^{-1}(x)$ are graded left $H_{[*]}(Z)$-module. \\
   We call these the \textbf{Springer fibre modules}.
\end{itemize}
Similarly in all examples one can obtain right module structure (the easy swaps are left to the reader). Independently, one can define the same graded module structure on $H_{*}(\pi^{-1}(x)), H^*(\pi^{-1}(x))$ using the description of the Steinberg algebra as Ext-algebra and a Yoneda operation (for this see
\cite{CG}, 8.6.13, p.448 ). 

\noindent
There is also a result of Joshua (see \cite{Jos}) saying that all hypercohomology groups $\HH_A^*(Z,F^\bullet), F^\bullet\in D_A^b(Z)$ carry the structure of a left (and right) 
$H_*^A(Z)$-module.

%%%%%%%%%%%%%%%%%%%%%%%%%%%%%%%%%%%%%%%%%%%%%%%%%%%%%%%%%%%%%%%%%%%%%%%%%%%%%%%%%%%%%%%%%
\subsection*{The Steinberg algebra}%, $H_A^*(pt)$ and $H_A^*(E)$}

\subsubsection*{The Steinberg algebra $H_{[*]}^A(Z)$ as module over $H^{-*}_A(pt)$.}

We set $\widetilde{W}:=\bigsqcup_{i,j\in I} W_{i,j}$ with $W_{i,j}:=W_i \setminus W /W_j$ where $W$ is the Weyl group for $(G,T)$ and $W_i\subset W$ is the Weyl group 
for $(L_i,T)$ with $L_i\subset P_i$ is the Levi subgoup. We will fix representatives $w\in G$ for all elements $w\in\widetilde{W}$. \\
Let $C_w= G\cdot (eP_i, wP_j)$ be the $G$-orbit in $G/P_i \times G/P_i$ corresponding to $w\in W_{i,j}$. 
\begin{lemma}
\begin{itemize}
\item[(1)] 
$p\colon C_w \subset G/P_i \times G/P_j \xrightarrow{pr_1} G/P_i$ is $G$-equivariant, locally trivial with 
fibre $p^{-1}(eP_i)= P_iwP_j/P_j$. 
\item[(2)]
$P_iwP_j/P_j$ admits a cell decomposition into affine spaces via Schubert cells $xB_jx^{-1}vP_j/P_j,  v\in W_i$ (and for a fixed $x\in W$ such that ${}^xP_j=P_i$, $B_j\subset P_j$ the Borel subgroup). In particular, $H_{odd}(P_iwP_j/P_j)=0$ and 
\[H_*(P_iwP_j/P_j)=\bigoplus_{v\in W_i}\C b_{i,j}(v), \quad  b_{i,j}(v):= [xB_jx^{-1}vP_j/P_j]. \]
It holds $\deg b_{i,j}(v)=2\ell_{i,j}(v)$ where $\ell_{i,j}(v)$ is the length of a minimal coset representative in $W$ for $x^{-1}vW_j \in W/W_j$.  
\item[(3)] For $A\in \{pt, T, G\}$ it holds $H_{odd}^A(C_w)=0$ and  since $G/P_i$ is simply connected 
\[H_n^A(C_w)=\bigoplus_{p+q=n}H^p_A(G/P_i) \otimes H_q(P_iwP_j/P_i), \quad  H_*^A(C_w)= \bigoplus_{u\in W/W_i, v\in W_i} \C b_{i}(u)\otimes b_{i,j}(v),\]
where $b_i(u)=[B_iuP_i/P_i]^*$ is of degree $2\dim_{\C}G/P_i-2\ell_i(u)$ with $\ell_i(u)$ is the length of a minimal coset representative for $u\in W/W_i$ and $b_{i,j}(v)$ as in (2). 
\end{itemize}
\end{lemma}

\paragraph{proof:} Is left out because it is standard techniques. %For sake of completeness we give a proof in the Appendix, see...

This implies the following properties for the homology of $Z$. 
\begin{lemma} \label{oddVan}
\begin{itemize}
\item[(1)] $Z$ has a filtration by closed $G$-invariant subvarieties such that the successive complements are 
$Z_w:= m^{-1}(C_w), w\in \widetilde{W}$ and the restriction of $m$ to $Z_w$ is a vector bundle over $C_w$ of rank $d_w$ (as complex vector bundle). 
Furthermore, 
\[
\begin{aligned}
 H_n^A(Z) &=\bigoplus_{w\in \widetilde{W}} H_n^A(Z_w) =\bigoplus_{w\in \widetilde{W}} H_{n-2d_w}^A (C_w)\\  
&=  \bigoplus_{i,j\in I}\bigoplus_{w\in W_{i,j}} \bigoplus_{u\in W/W_i, v\in W_i}\C b_i(u)\otimes b_{i,j}(v)  
\end{aligned}
\]
where the last direct sum goes over the $u,v$ with the property $2\dim G/P_i -2\ell_i(u) +2\ell_{i,j}(v) = n-2d_w$.
\item[(2)] $H_{odd}(Z)=0, H^{odd}(Z)=0.$
\item[(3)] $Z$ is equivariantly formal (for $T$ and $G$, for Borel-Moore homology and cohomology).\\
In particular, for $A\in \{T,G\}$ the following forgetful maps are 
surjective $H_*^A(Z)\surj H_*(Z)$ and $H^*_A(Z)\surj H^*(Z)$
algebra homomorphisms. It even holds the stronger isomorphism of $\C$-algebras 
\[
\begin{aligned}
H_*(Z)    %= H_*^G(Z)\otimes_{H_*^G(pt)} (H_*^G(pt)/H_{>0}^G(pt))  \\
            & = H_*^A(Z) / H_{<0}^A(pt)H_*^A(Z)\\
            %& = H_*^T(Z) \otimes_{H_*^T(pt)} (H_*^T(pt)/H_{>0}^T(pt) \\
            %& = H_*^T(Z)/ H_{<0}^T(pt)H_*^T(Z) \\
H^*(Z)   %& = H^*_G(Z)\otimes_{H^*_G(pt)} (H^*_G(pt)/H^{>0}_G(pt))  
           &= H^*_A(Z) / H^{>0}_A(pt)H^*_A(Z)\\
         %& = H^*_T(Z) \otimes_{H^*_T(pt)} (H^*_T(pt)/H^{>0}_T(pt) 
           %&= H^*_T(Z)/ H^{>0}_T(pt)H^*_T(Z)                   
\end{aligned}
\]
As a consequence we get the following isomorphisms. 
\[
\begin{aligned}
1)\;  H_*^A(Z) & = H_*(Z)\otimes_{\C} H_*^A(pt)   \quad \text{ of }H_*^A(pt)\text{-modules} \\
%2)\;   H_*^T(Z) & = H_*(Z)\otimes_{\C} H_*^T(pt) \quad \text{ of }H_*^T(pt)\text{-modules}\\
2)\;  H^*_A(Z) & = H^*(Z)\otimes_{\C} H^*_A(pt) \quad \text{ of }H^*_A(pt)\text{-modules}\\
%4)\;  H^*_T(Z) & = H^*(Z)\otimes_{\C} H^*_T(pt) \quad \text{ of }H^*_T(pt)\text{-modules}
\end{aligned} 
\] 
We can see that $H_{[*]}^A(Z)$ has finite dimensional graded pieces and the graded pieces are bounded from below in negative degrees.
\end{itemize}
\end{lemma}

%%%%%%%%%%%%%%%%%%%%%%%%%%%%%%%%%%%%%%%%%%%%%%%%%%%%%%%%%%%%%%%%%%%%%%%%%%%%%%%%%%%%%%%%%%%%%%%%%%%%%%%%%%%%%%%%%%%%%%%%%%%%%%%%
\subsubsection*{The Steinberg algebra $H_*^A(Z)$ and $H_A^*(E)$ }

%It holds $Z\circ E = E$, therefore we have an operation $*$ of $H(Z)$ on $H(E)$. 
Recall from a previous section that $H_A^*(E)$ is a graded left (and right) $H_{[*]}^A(Z)$-module 
and that $H_A^*(E)$ has a $H_A^*(pt)$-algebra structure with respect to the cup product, the $H_{[*]}^A(Z)$-operation is $H_A^*(pt)$-linear. 
%%Then, the Poincare duality $H_A^*(E)\cong H_*^A(E)$ make $H_A^*(E)$ a module under $H_*^A(Z)$, this operation is $H_A^*(pt)$-linear.

\begin{rem}
Let $q_i\colon E_i\to pt, \; i\in I$, there is an isomorphism of algebras 
\[ \End_{H_A^*(pt)}(H^*_A(E)) = H_{*}^A(E\times E)  = \Ext^*_{D^A(pt)} (\bigoplus_{i\in I} {(q_i)}_*\underline{\C}[e_i], \bigoplus_{i\in I}{(q_i)}_*\underline{\C}[e_i]),\]
the first equality follows from \cite{CG},Ex. 2.7.43, p.123, 
for the second first use the Thom isomorphism to replace $E\times E$ by a union of flag varieties, then use 
theorem \ref{Extalgebra} from Chriss and Ginzburg for the Springer map given by the projection to a point. \\
Furthermore, under the identifications, the following three algebra homomorphisms are equal. 
\begin{itemize}
\item[(1)] The map $H_*^A(Z)\to \End_{H^*_A(pt)}(H^*_A(E)), c\mapsto (e\mapsto c*e)$. 
\item[(2)] $i_*\colon H_*^A(Z) \to H_*^A(E\times E)$ where $i\colon Z\to E\times E$ is the natural embedding.  
\item[(3)] $\begin{aligned}\Ext^*_{D^A(V)} (\bigoplus_{i\in I}{(\pi_i)}_*\underline{\C}[e_i], \bigoplus_{i\in I}{(\pi_i)}_*\underline{\C}[e_i]) &\to 
\Ext^*_{D^A(pt)} (a_*(\bigoplus_{i\in I}{(\pi_i)}_*\underline{\C})[e_i], a_*(\bigoplus_{i\in I}{(\pi_i)}_*\underline{\C}[e_i])),\\
 f&\mapsto a_*(f)\end{aligned}$\\
where $a\colon V\to pt$.  
\end{itemize}
\end{rem}

\begin{lemma}\label{inj} (\cite{VV2}, remark after Prop.3.1, p.12)
Assume that $T\subset \bigcap_i P_i$ is a maximal torus 
and $Z^T= E^T\times E^T, E^T= \bigsqcup_{i\in I}(G/P_i)^T$.
Let $A\in \{T,G\}$. There is an injective homomorphism of $H_A^*(pt)$-algebras 
\[ H_*^A(Z) \to \End_{H_A^*(pt)} (H_A^*(E)),\] 
Let $\mathfrak{t}$ be the Lie algebra of $T$, then it holds $H_G^*(E)\cong\C[\mathfrak{t}]^{\oplus I}$, where $\C[\mathfrak{t}]$ is the ring of regular functions on the affine space $\mathfrak{t}$.  
%which is implied by the Goresky-Kottwitz-MacPherson Theorem. 
\end{lemma}

\paragraph{proof:}
For $G$-equivariant Borel-Moore homology we claim that the following diagram is commutative 
\[
\xymatrix{
H_*^T(Z^T) \ar[r] & H_*^T(E^T\times E^T) \\
H_*^T(Z) \ar[r]\ar[u] & H_*^T(E\times E) \ar[u] \\
H_*^G(Z) \ar[r] \ar[u] & H_*^G(E\times E) \ar[u]
}
\]
The commutativity of the lowest square uses functoriality of the forgetful maps. By assumption $Z^T=(E\times E)^T$, the highest horizontal map is an 
isomorphsim. Now, by the GKM-localization theorem the two vertical maps in the upper square are injective. That implies that the middle horizontal map 
has to be injective, together with (2) from the previous remark it implies the claim for $T$-equivariant Borel-Moore homology. 
But by the splitting principle, i.e. the identification of the $G$-equivariant Borel-Moore homology with the $W$-invarient subspace in the $T$-equivariant Borel-Moore homology, the forgetful maps become the inclusion of the $W$-invariant subspace. This means the two vertical maps in the lower square are injective. 
This implies that the lowest horizontal map is injective. Together, with (2) of the previous remark the claim follows for $A=G$. \hfill $\Box$

The main ingredient to the previous lemma is Goretzky's,Kottwitz' and MacPherson's localization theorem (see \cite{GKM}). 
Similar methods are currently developped by Gonzales for $K$-theory in \cite{Go}. 

The previous lemma is wrong for not equivariant Borel-Moore homology as the following example shows. 
\begin{exa}
Let $G$ be a reductive group with a Borel subgroup $B$ and $\mathfrak{u}$ be the Lie algebra of its unipotent radical. 
$Z:= (G\times^B\mathfrak{u})\times_{\mathfrak{g}} (G\times^B\mathfrak{u})$, then it holds that the algebra $H_*(Z)$ can under the isomorphism 
in Kwon (see \cite{Kw}) be identified with $\C[\mathfrak{t}]/I_W \# \C[W]  $ which is defined as the $\C$-vector space $\C[\mathfrak{t}]/I_W\otimes_{\C} \C[W]$  
with the multiplication $(f\otimes w)\cdot (g\otimes v) := fw(g)\otimes wv$. Furthermore, we can identify $\End_{\C}(H^*(E))$ via the Thom-isomorphism and 
the Borel map with $\End_{\C -lin}(\C[\mathfrak{t}]/I_W)$. The canonical map identifies with 
\[
\begin{aligned}
\C[\mathfrak{t}]/I_W \# \C[W] & \to  \End_{\C -\lin}(\C[\mathfrak{t}]/I_W)\\
  f\otimes w &\mapsto (p\mapsto f w(p))
\end{aligned}
\]
This map is neither injective nor surjective. 
For example $\sum_{w\in W} 1\otimes w \neq 0 $ in $\C[\mathfrak{t}]/I_W \# \C[W]$ but its image 
$(p\mapsto \sum_{w\in W} w(p))$ is zero because $\sum_{w\in W} w(p) \in I_W$. Because both spaces have the same $\C$-vector space dimension, it is clear 
that it is also not surjective. 
\end{exa}

%%%%%%%%%%%%%%%%%%%%%%%%%%%%%%%%%%-H(Z) as H(E)-module

Furthermore, $H_*^A(Z)$ is naturally a $H_A^*(E)$-module. It holds 
In fact, $H^*_A(E)\cong H_*^A(\bigsqcup_{e_{ij}=\overline{e}\in W_{i,j}} Z^{e_{i,j}})$ is even a subalgebra of $H_*^A(Z)$.

\begin{coro} \label{centre} In the situation of the previous lemma, i.e. $T\subset \bigcap_i P_i$ is a maximal torus 
and $Z^T= E^T\times E^T, E^T= \bigsqcup_{i\in I}(G/P_i)^T$ and let $A\in \{T,G\}$. 
There are injective homomorphism of $H_G^*(pt)$-algebras 
\[ H_*^A(pt) \subset H_*^A(E) \subset  H_*^A(Z) \to \End_{H_A^*(pt)} (H_A^*(E)),\] 
where the first inclusion is given by the pullback along the map $E\to pt$. In particular, $H_A^*(pt)$ is contained in the centre of $H_*^A(Z)$ (we only know examples where it is equal to the centre). 
\end{coro}
 
Let $w\in \widetilde{W}$. 
Observe, that $H_A^*(E)$ already operates on $H_*^A(Z^w)$ and the composition 
$ H_*^A(Z)= \bigoplus_{w} H_*^A(Z^w)$ is a direct sum composition of $H_A^*(E)$-modules. 
Using the Thom-isomorphism, up to a degree shift we can also study $H_*^A(C^w)$ as 
module over $H^*_A(\bigsqcup_{I} G/P_i)$. Now, let $e_i$ be the idempotent in $H^*_A(\bigsqcup_{I} G/P_i)=\bigoplus_{i\in I}H^*_A(G/P_i)$ which corresponds to 
the projection on the $i$-th direct summand. Since for $w\in W_{i,j}$ it holds $H_*^A(C_w)=H^*_A(G/P_i)\otimes_{\C} H_*(P_iwP_j/P_j)$ also as $H_A^*(G/P_i)$-module, we conclude that $H_*^A(C_w)$ is always a projective module over $H_A^*(\bigsqcup_{I} G/P_i)$.

\begin{lemma}
\begin{itemize}
\item[(1)] Let $w\in W_{i,j}$. 
Each $H_*^A(Z^w)$ is a projective graded $H^*_A(E)$-module of the form 
\[\bigoplus_{ v\in W_i} (H_A^*(E)e_i) [2d_w + \deg b_{i,j}(v)],\]
 where $[d]$ denotes the degree shift by $d$. In particular, $H_*^A(Z)$ is a projective graded $H^*_A(E)$-module. 
\item[(2)] If all $P_i=B_i$ are Borel subgroups of $G$, then $H_*^A(Z)= \bigoplus_{w,j \in W\times I}(\bigoplus_{i\in I} (H_A^*(E)e_i)[d_{w,i,j}])$ as graded $H_A^*(E)$-module for certain $d_{w,i,j} \in \Z$. In particular, if we forget the grading $H_*^A(Z)$ is a free $H_A^*(E)$-module of rank 
$\# W\cdot\# I$. 
\end{itemize}
\end{lemma}

%\paragraph{proof:}

%\hfill $\Box $
%It also holds $H_*(Z) = H_*^T(Z)/ mH_*^T(pt)$ with $m=\Ker ..$ and $H_*^G(Z)=(H_*^T(Z))^W$ as convolution algebras. 

%%%%%%%%%%%%%%%%%%%%%%%%%%%%%%%%%%%%%%%%%%%%%%%%%%%%%%%%%%%%%%%%%%%%%%%%%%%%%%%%%%%%%%%%%%%%%%%%%%%%%%%%%%%%%%%%%%%%%%%%%%%%%%%%%%%%%%%%%%%%%%%%

\subsection*{Indecomposable projective graded modules over $H_{[*]}^A(Z)$ and their tops for a different grading}
%\subsubsection{Modules from the decomposition theorem}
%Here, the grading is the Ext-grading on $H_*^A(Z)$ from remark \ref{grading}.

Let $X$ be an irreducible algebraic variety, 
we call a decomposition $X=\bigcup_{a\in \mcA} S_a$ into finitely many irreducible smooth locally closed subsets a \emph{weak} stratification.  
Since $\pi\colon E=\bigsqcup_{i\in I} E_i \to V$ is a $G$-equivariant projective map, there exists (and we fix it) a \emph{weak} stratification into $G$-invariant subsets $V=\bigsqcup_{a\in \mcA}{S_a}$ such that $\pi^{-1}(S_a) \xrightarrow{\pi} S_a$ is a locally trivial\footnote{with respect to the analytic topology} fibration with constant fibre 
$F_a:= \pi^{-1}(s_a)$ where $s_a\in S_a$ one fixed point, for every $a\in \mcA$. (For projective maps of complex algebraic varieties one can always find such a weak stratification, see \cite{Ar3}, 4.4.1-4.4.3)\footnote{If the image of $\pi$ is irreducible, 
by \cite{Ar3}, theorem 1.9.10 we can refine this stratification to a (finite) Whithney stratification, but it is not clear if we can find a Whitney stratification into $G$-invariant subsets.}\\

Recall that for any $G$-equivariant projective map of complex varieties, 
the decomposition theorem (compare \cite{BBD} for the not equivariant version and \cite{BL} for the equivariant version) 
holds. Let $A\in \{ pt, T,G\}$, we denote by $D_A^b(V)$ the $A$-equivariant derived category defined by Bernstein and Lunts in \cite{BL}. Let $t$ run over all irreducibel $G$-equivariant local systems $\mcL_t$ on some stratum $S_t=S_{a_t}$, $a_t\in \mcA$, 
we write $IC_t^A:=(i_{\overline{S_t}})_*(\mathcal{I}\mathcal{C}^A(S_t,\mcL_{t})[d_{S_t}]$ with $d_{S_t}=\dim_{\C-var} S_t$  for the simple perverse sheaf in the category of $A$-equivariant perverse sheaves $Perv_A(V) \subset D_A^b(V)$, see again \cite{BL}, p. 41. 
Let $e_i=\dim_{\C} E_i, i\in I$, then $ \underline{\C}_{E_i}[e_i]$ is a simple perverse sheaf in $D_A^b(E)$.
For a graded vector space $L=\bigoplus_{d\in Z}L_d$ we define $L(n)$ to be the graded vector space with $L(n)_d:= L_{n+d}$, we see $\C$ as the graded 
vector space concentrated in degree zero, $n\in \Z$.   
For an element $F^{\bullet} \in D_A^b(X)$ for an $A$-variety $X$ we write $F^{\bullet}[n]$ for the (class of the) complex $(F^{\bullet}[n])_d:= F^{d+n}$, $n\in \Z$.   
Now given $F^{\bullet}\in D_A^b(X)$ 
and a finite dimensional graded vector space $L:= \bigoplus_{i=1}^r \C(d_i)$ we define 
\[ L \otimes_{gr} F^{\bullet} := \bigoplus_{i=1}^r F^{\bullet}[d_i] \in D_A^b(X) \]  
The $A$-equivariant decomposition theorem applied to $\pi$ gives
%applied to $\pi$ gives a direct sum decomposition in $D_A^b(V)$ (in the sense of \cite{BL}\footnote{in fact the decomposition takes place in the full 
%subcategory of complexes of sheaves on $X_A$ which lie in $D_A^b(X)$ and whose cohomology is locally constant over the strata $(S_a)_A, \; a\in A$})  
\[
\begin{aligned}
 \bigoplus_{i\in I} (\pi_i)_*\underline{\C}_{E_i}[e_i] &= \bigoplus_{t} L_{t}\otimes_{gr} IC_t^A \quad \in D^b_A(V)  
%\pi_*(\underline{\C}[d]) &= \bigoplus_{k\in \Z,\; t} L_{t}(k)\otimes_{\C}IC_t^T[k] \quad \in D^b_T(V)  \\
%\pi_*(\underline{\C}[d]) &= \bigoplus_{k\in \Z,\; t} L_{t}(k)\otimes_{\C}IC_t^G[k] \quad \in D^b_G(V)  \\
\end{aligned}
\]
where the $L_t:= \bigoplus_{d\in \Z}L_{t,d}$ are complex finite dimensional graded vector spaces.\\
Let $\D$ be the Verdier-duality on $V$, it holds $\D (\pi_*(\underline{\C}[d]))= \pi_*(\underline{\C}[d]), \D (IC_t^A) = IC_{t^*}^A$ where for $t=(S, \mcL)$ it 
holds $t^*=(S, \mcL^*),\;\; \mcL^*:=\mcH om (\mcL, \underline{\C} )$. This implies $L_t=L_{t^*}$ for all $t$.

\subsubsection*{Indecomposable projectives in the category of graded left $H_{[*]}^A(Z)$-modules}

%\noindent
%We will only consider categories of finitely generated graded modules where the grading is the Ext-grading on $H_*^A(Z)$. 
We set
\[
P_{t}^A:= \Ext_{D_A^b(V)}^*(IC_t^A,\bigoplus_{i\in I}{(\pi_i)}_* \underline{\C}[e_i]).
\] 
It is a graded (left) $H_{[*]}^A(Z)$-module. It is imdecomposable because $IC_t^A$ is simple.  
Clearly it holds as left graded $H_{[*]}^A(Z)$-modules  
\[
\begin{aligned}
H_{[*]}^A(Z) %&=  \bigopus_{n\in \Z} \bigoplus_{m\in \Z, \; s } Ext_{D_A^b(V)}^n(\pi_*\underline{\C}, L_s(m) \otimes IC_s^A[m])  \\
         &= \bigoplus_{d \in \Z, t} L_{t,d} \otimes [\bigoplus_{n\in \Z}  \Ext_{D_A^b(V)}^{n+d}(IC_t^A, \bigoplus_{i\in I}{(\pi_i)}_* \underline{\C}[e_i])]\\
         &= \bigoplus_{d \in \Z, t}  L_{t,d} \otimes_{\C} P_t^A[d]\\
         %&= \bigoplus_{m,n\in \Z , s} [ L_s(m) \otimes   Ext_{D_A^b(V)}^{n}( .., IC_s^A) ]\\
&= \bigoplus_{t} L_t\otimes_{gr} P_{t}^A
\end{aligned}
\]
that implies that $P_t^A$ is a projective module and that $(P_t^A)_t$ is a complete set of isomorphism 
classes up to shift of indecomposable projective graded $H_{[*]}^A(Z)$-modules. 

Assume that $H_A^*(pt)$ is a graded subalgebra of the centre of $H_{[*]}^A(Z)$, compare corollary \ref{centre}. 
\begin{lemma}
The elements $H_A^{>0}(pt)$ operate on any graded simple $H_{[*]}^A(Z)$-module $S$ by zero. 
In particular, by lemma \ref{oddVan} we see that $S$ is a graded simple modules over $H_{[*]}(Z)$. Any graded simple module 
is finite-dimensional and there exists up to isomorphism and shift only finitely many graded simple modules. \\
For any graded simple module $S$ there is no nonzero degree zero homomorphism $S\to S(a), \;a\neq 0$.   
\end{lemma}

\paragraph{proof:}
Since $S$ is simple it holds $H_A^{>0}(pt)\cdot S$  is zero or $S$. Assume it is $S$, pick a non-nilpotent element $x\in H_A^d(pt)$ and $y\in S,y\neq 0$, homogeneous. Then, it holds $S= H_{[*]}^A(Z)\cdot y= H_{[*]}^A(Z)\cdot x^n y$ contradicting the fact that $S$ has to have a minimal degree generator. 
Therefore $H_A^{>0}(pt)\cdot S=0.$\\
By \cite{NvO}, II.6, p.106, we know that the graded simple modules considered as modules over the ungraded rings $H_*^A(Z), H_*(Z)$ are still simple modules. 
Since the finite-dimensional algebra $H_*(Z)$ has up to isomorphism only finitetly many simples, the claim follows.   \\ 
Any nonzero degree $0$ homomorphism $\phi\colon S\to S(a)$ has to be an isomorphism. Let $S=H_{[*]}^A(Z)\cdot y$ as before, set $\deg y=m$. Then $S(a)= H_{[*]}^A(Z)\cdot \phi(y)$, $\deg\; \phi (y)=m$ which gives a contradictions when considering the minimal nonzero degrees of $S$ and $S(a)$.
\hfill $\Box $ 

\begin{coro} \label{projectives}
There is a bijection between isomorphism classes up to shift of  
\begin{itemize}
\item[(1)] indecomposable projective graded $H_{[*]}^A(Z)$-modules
\item[(2)] indecomposable projective graded $H_{[*]}(Z)$-modules
\item[(3)] simple graded $H_{[*]}(Z)$-modules.
\end{itemize}
The bijection between (1) and (2) is clear from the decomposition theorem, it maps $P\mapsto P/H_A^{>0}(pt)P$. 
We pass from (3) to (2) by taking the projective cover and we pass 
from (2) to (3) by taking the top 
(which is graded because for a finite dimensional graded algebra the radical is given by a graded ideal). 
\end{coro}

\begin{exa} (due to Khovanov and Lauda, \cite{KL})
Let $G\supset B\supset T$ be a reductive group containing a Borel subgroup containing a maximal torus, $Z= G/B \times G/B$. 
Then, it is known that $H_*^G(Z)= \End_{\C[\mathfrak{t}]^W }(\C[\mathfrak{t}])=:NH$ where $W$ is the Weyl group associated to $(G,T)$ and 
$\mathfrak{t}=Lie(T)$. \\
The $G$-equivariant pushforward (to the point) of the shift of the constant sheaf is a direct sum of shifts of copies of the constant sheaves on the point, therefore there exist precisely one indecomposable projective graded $H_{[*]}^G(Z)$-module up to isomorphism and shift. It is easy to see that $P:=\C[\mathfrak{t}]$ is an indecomposable projective module and $P/H_G^{>0}(pt)P=\C[\mathfrak{t}]/I_W$ is the only graded simple $NH$-module which is the top of $P$. 
Also, one checks that $H_{[*]}(Z)=\End_{\C}(\C[\mathfrak{t}]/I_W)$ is a semi-simple algebra which has up to isomorphism and shift only the one graded 
simple module $\C[\mathfrak{t}]/I_W$.    
\end{exa}

Now we equipp the Steinberg algebra with a grading by positive integers which leads to a description of graded simple modules in terms of the multiplicity vector 
spaces $L_t$ in the BBD-decomposition theorem. 

\subsubsection*{Simple objects in the category of graded finitely generated left $H_{<*>}^A(Z)$-modules}

Given a graded vector space $L$, we write $\langle L\rangle :=\bigoplus_{d\in \Z}L_d$ for the underlying (ungraded) vector space.  
%This is a citation from \cite{CG}, thm 8.6.12, p.448 and the proof from page 447. \\
If we regrade $H_*^A(Z)$ as follows 
\[ H_{<n>}^A(Z) :=  \bigoplus_{s,t} \Hom_{\C}(\langle L_t \rangle , \langle L_s \rangle ) \otimes \Ext_{D_A^b(V)}^{n}(IC_t^A, IC_s^A)], \]
in other words 
\[H_{<*>}^A(Z)=\Ext^*(\bigoplus_{t} \langle L_{t}\rangle \otimes_{\C} IC_t^A, \bigoplus_{t} \langle L_{t}\rangle  \otimes_{\C}IC_t^A)\]
as graded algebra. It holds as an algebra this one is isomorphic to $H_*^A(Z)$. With the same arguments as in the previous 
section one sees that $P_t^A:= \Ext_{D_A^b(V)}^*(IC_t^A,\pi_* \underline{\C})$ are a complete representative system for the isomorphism classes of the indecomposable projective graded $H_{<*>}^A(Z)$-modules.  
 
\noindent
We claim that there is a graded $H_{<*>}^A(Z)$-module structure on the (multiplicity-)vector space $\langle L_t\rangle $ such 
that the family $\{\langle L_t\rangle \}_t$ is a complete set of 
the isomorphism classes up to shift of graded simple modules. 
Using $\Hom (IC_t^A, IC_s^A) =\C \delta_{s,t}, \; \Ext^n(IC_t^A, IC_s^A)=0$ for $n<0$ we get 
%%To see this consider  
%\[
%\begin{aligned}
%H_*^A(Z) &=  \bigopus_{n\in \Z} \bigoplus_{d,m\in \Z, \; t,s } Ext_{D_A^b(V)}^n((L_t)_d\otimes IC_t^A[k], L_s(m) \otimes IC_s^A[m])  \\
%%         &= \bigoplus_{k,m \in \Z, s,t} \Hom_{\C}(L_t(k), L_s(m)) \otimes [\bigoplus_{n\in \Z}  Ext_{D_A^b(V)}^{n-m+k}(IC_t^A, IC_s^A)]\\
%%        % &= \bigoplus_{k,m \in \Z, s,t} \Hom_{\C}(L_t(k), L_s(m)) \otimes \bigoplus_{n\in \Z}  Ext_{D_A^b(V)}^{n}(IC_t^A, IC_s^A)[k-m]\\
%%         &= \bigoplus_{s,t} [\Hom (L_t, L_s) \otimes   Ext_{D_A^b(V)}^{*}(IC_t^A, IC_s^A)]
%%\end{aligned}
%\]
%where we see $\Hom (L_t, L_s)= \bigoplus_{k,m\in \Z} \Hom_{\C}(L_t(k), L_s(m))$ as a graded vector space with \\
%$\deg  \Hom_{\C}(L_t(k), L_s(m)) = k-m$ and 
%the last tensor product is a tensor product of graded spaces. 
\[
H_{<*>}^A(Z) = \underbrace{\bigoplus_t \End (\langle L_t \rangle )}_{\deg =0} \oplus \bigoplus_{s,t} \Hom (\langle L_t\rangle ,\langle L_s\rangle ) \otimes_{\C} 
\Ext^{>0} (IC_t^A, IC_s^A).
\]  
Now, the \emph{second} summand is the graded radical, 
i.e. the elements of degree $>0$ (with respect to the new grading). It follows 
\[
H_{<*>}^A(Z) \surj H_{<*>}^A(Z)/(H_{<*>}^A(Z))_{>0} = \bigoplus_{t}\End_{\C} (\langle L_t\rangle ).
\]
This gives $\langle L_t\rangle $ a natural graded $H_{<*>}^A(Z)$-modul structure concentrated in degree zero (the positive degree elements in $H_{<*>}^A(Z)$ operate by zero). 
Observe, that $\langle L_t\rangle $ does not depend on $A$, i.e. in fact they are modules over $H_{<*>}^A(Z)$ via the forgetful morphism $H_{<*>}^A(Z) \to H_{<*>}(Z)$. \\
That means we can instead look for the simple graded modules of $H_{<*>}(Z)$. 

\begin{rem} Let $H_*$ is a finite dimensional positively graded algebra such that  
\[ H_0= H_*/H_{>0} =\bigoplus_{t}\End (L_t) \]
is a semi-simple algebra. Then $H_{>0}$ is the set of nilpotent elements, i.e. Jacobsen radical of $H_*$. 
Furthermore all simple and projective $H_*$-modules are graded modules. 
\begin{itemize}
\item[*] $(L_t)_t$ is the tuple of (pairwise distinct isomorphism classes of all) simple modules. 
\item[*]
For each $t$ pick an $e_t\in \End (L_t) \subset H_0$ which corresponds to projection and then inclusion of a one dimensional subspace of $L_t$. \\
$(P_t:= H_* \cdot e_t)_t$ is the tuple of (pairwise distinct isomorphism classes of all) indecomposable projective modules.  
\end{itemize}
\end{rem}
We can apply this remark to $H=H_{<*>}(Z)$. As a consequence we see that up to shift $(\langle L_t \rangle )_t$ is the tuple of (pairwise distinct isomorphism classes of all) simple graded $H_*^A(Z)$-modules.

From now on, the case where the two gradings coincide will play a special role. 
\begin{rem} The following conditions are equivalent 
\begin{itemize}
\item[(1)] $H_{[*]}^A(Z)=H_{<*>}^A(Z)$ as graded algebra for every (at least one) $A\in \{pt, T, G\}$. 
\item[(2)] $(\pi_i)_*\underline{\C}[e_i]$ is $A$-equivariant perverse for every $i\in I$ for every (at least one) $A\in \{pt, T, G\}$. 
\item[(3)] $\pi_i\colon E_i \to V$ is \textbf{semi-small} for every $i\in I$, this means by definition $\dim Z_{i,i}=e_i$ for every $i\in I$. 
\end{itemize}
In this case, we say the Springer map is semi-small. Also, $\pi$ semi-small is equivalent to $H_{top}(Z_{i,i})=H_{[0]}(Z_{i,i}), i\in I$. Observe, that $H_{[0]}(Z)$ is always a subalgebra and in the semi-small case isomorphic to the quotient algebra 
$H_{[*]}(Z)/  (H_{[*]}(Z))_{>0}$. 
Assume $\pi$ semi-small, then it holds $2\dim \pi_i^{-1}(x)\leq e_i -d_S$ where $x\in S$ belongs to the stratification and 
$H_{top}(\pi^{-1}(x)):= \bigoplus_{i\colon 2\dim \pi_i^{-1}(x) = e_i - d_{S}}H_{2\dim \pi_i^{-1}(x)}(\pi_i^{-1}(x))$ is a left $H_{[0]}(Z)$-modules via the restriction of the 
convolution construction. If $I$ consists of a single element, $H_{top}(Z) = H_{[0]}(Z)$ and $H_{2\dim \pi^{-1}(x)}(\pi^{-1}(x))$ is an $H_{[0]}(Z)$-module independent of the 
condition $2\dim \pi^{-1}(x) = \dim E - d_{S}$.
\end{rem}

\begin{rem}
If one applies the decomposition theorem to $\pi_i, i\in I$ one gets that $L_t=\bigoplus_{i\in I} L_t^{(i)}$ (as graded vector space) where 
$L_t^{(i)}$ is the multiplicity vector space for $IC_t$ in $(\pi_i)_*\underline{\C}[e_i]$. It holds $\{L_t^{(i)}\mid L_t^{(i)}\neq 0\}$ is the complete set of isomorphism classes of simple $H_*(Z_{i,i})$-modules. 
\end{rem}

\begin{rem} In fact, Syu Kato pointed out that the categories of finitely generated graded modules over $H_{[*]}^A(Z)$ and $H_{<*>}^A(Z)$ are equivalent. 
This has been used in \cite{Ka3}. 
\end{rem}

\begin{rem}
Now, we know that the forgetful (=forgetting the grading) functor from finite dimensional graded $H_{[*]}(Z)$-modules to finite dimensional 
$H_*(Z)$-modules maps graded simple modules to simple modules. We can use the fact that we know that simples and graded simples are parametrized by 
the same set to see: Every simple $H_*(Z)$-module $L_t$ has a grading such that it becomes a graded simple $H_{[*]}(Z)$-module and every graded simple is of this form.   
\end{rem}

For the decomposition matrix for the finite dimensional algebra $H_*(Z)$, there is the following result of Chriss and Ginzburg. 
\begin{satz} (\cite{CG}, thm 8.7.5) 
Assume $H_{odd}(\pi^{-1}(x))=0$ for all $x\in V$. Then, the following matrix multiplication holds 
\[ [P\colon L] = IC \cdot D\cdot IC^t \]
where all are matrices indexed by $s=(S,\mcL ), t =(S^\prime, \mcL^\prime )$ such that $L_t\neq 0, L_s\neq 0$ 
and $()^t$ denote the transposed matrix.  
\[
\begin{aligned}
{}[ P \colon L ]_{s,t} &:= [P_s\colon L_t] = \sum_k \dim \Ext^k(IC_t, IC_s)\\
IC_{s,t}          &:= \sum_{k} [\mcH^k (i_S^*(IC_t)) \colon \mcL ]                   \\
D_{s,t}           &:= \delta_{S,S^\prime} \sum_k (-1)^k \dim H^k(S, (\mcL^\prime)^* \otimes \mcL ) 
\end{aligned}
\]
\end{satz} 

\noindent

According to Kato in \cite{Ka3}, the whole theory of these algebras is reminiscent of quasi-hereditary algebras (but we have infinite dimensional algebras). 
He introduces standard and costandard modules for $H_{<*>}^G(Z)$ in \cite{Ka3},thm 1.3, under some assumptions.% \footnote{In this version of \cite{Ka3}, it becomes not clear that the grading $H_{<*>}^G(Z)$ is not equal to the grading $H_{[*]}^G(Z)$, see proof of lemma 1.2. 
% As far as I understand it he uses the fact that the two algebras are Morita 
% equivalent to deduce that the categories of graded modules are equivalent - which I doubt is true.}
He shows that under these assumptions\footnote{ $ = $ finitely many orbits with connected stabilizer groups in the image of the Springer map, \emph{pure of weight zero} for $H_{<*>}^G(Z)$ and of the $IC_t$ in the decomposition theorem}, $H_{<*>}^G(Z)$ has finite global dimension (see \cite{Ka3}, thm 3.5).

%Other comments are that the theory is reminiscent of affine celluar algebras, of highest weight theories and Kazhdan-Luzstig theory but it is none of it really fit. 
% 
%When one considers the new grading there are also interesting modules called standard and costandard modules, and a lot of extra structure which is 
%reminiscent of the theory of . Some cases have been explored by S. Kato in \cite{Ka3} and by Brundan, Kleshchev, McNamara in \cite{BKM}. Chriss and Ginzburg saw the following theorem as an indication for that which gives a description of the decomposition matrix (i.e. the multiplicity $[P_s\colon L_t]$ of a simple $H_*(Z)$-module $L_t$ in a decomposition series of a indecomposable projective module $P_s$). 

\subsubsection*{Springer fibre modules in the category of graded $H_{?}^A(Z)$-modules}

Recall, that Springer fibre modules $H_{[*]}(\pi^{-1}(x)), H^{[*]}(\pi^{-1}(x)), x\in V$ are naturally graded modules over $H_{[*]}(Z)$, but if we forget about the grading and we can show 
that they are actually semi-simple or projective in $H_*(Z)$-mod, then, we can see them as semi-simple graded $H_{<*>}^A(Z)$-modules for $A\in \{G,T,pt\}$ by the previous section. 

Let $A=pt$. 
Since the map $\pi$ is locally trivial over $S:=S_a$ we find that 
\[ i_S^*( \bigoplus_{i\in I}R^k(\pi_i)_* \underline{\C}[e_i]),\quad i_S^{!}(\bigoplus_{i\in I}R^k(\pi_i)_* \underline{\C}[e_i] \] 
are local systems on $S$, via monodromy they correspond to the $\pi_1(S,s)$-representations
\[ H^{[k]}(\pi^{-1}(s))=\bigoplus_{i\in I}H^{e_i+k}(\pi_i^{-1}(s)) , \quad  \bigoplus_{i\in I}H_{e_i-k}(\pi_i^{-1}(s))=H_{[k]} (\pi^{-1}(s)) \]
with $e_i:= \dim_{\C} E_i$ respectively (for a fixed point $s=s_a\in S$, cp. \cite{CG}, Lemma 8.5.4).   \\
Now, let us make the extra assumption that the image of the Springer map is irreducible 
and the stratification $\{S_a\}_{a\in \mcA}$ is a Whitney stratification (every algebraic stratification of an irreducible variety can be refined to a 
Whitney stratification see \cite{Ar3}, thm 1.9.10, p.30), which is totally ordered by inclusion into the closure. Let $S\subset \overline{S^\prime}$ for two strata 
$S, S^\prime$, we write $\Ind_{S^\prime}^S(\mcL):= i_S^*\circ \mcH^* (IC_{(S^\prime , \mcL)}) $, i.e. we consider the 
functors for  $k\in [-d_{S^\prime}, -d_S]$  
\[
\begin{aligned}
\Ind_{S^\prime}^S (-)_k \colon LocSys(S^\prime) &\to LocSys(S)\\
                                       \mcL &\to \Ind_{S^\prime}^S (\mcL)_k:= i_S^*\circ \mcH^k(IC_{(S^\prime , \mcL)})
\end{aligned}
\]  
where $LocSys(S)$ is the category of local systems on $S$, i.e. locally constant sheaves on $S$ of finite dimensional vector spaces.
(for other $k\in \Z$ this is the zero functor). If we apply the functor $i_S^*\circ \mcH^k$ on the right hand side of the decomposition theorem 
we notice the following (for the cohomology groups of IC-sheaves, see \cite{Ar3}, section 4.1, p.41), let $t=(S^\prime , \mcL)$. 
\[  
i_S^*\mcH^k( IC_t)= \begin{cases}
%0,                               & \text{ if }  d_S  > d_{S^\prime} \\
%                                 & \text{ or }  d_S =  d_{S^\prime}  , i\neq d_S \\
%                                 & \text{ or }  d_S  < d_{S^\prime},   i\in (-\infty ,-d_{S^\prime}-1]\cup [-d_S, \infty)    \\
\mcL ,                       & \text{ if }  d_S =  d_{S^\prime} ,  k=-d_S \\
\Ind_{S^\prime}^S (\mcL)_k  & \text{ if }  d_S  < d_{S^\prime},   k\in [-d_{S^\prime}, -d_S-1]  \\
0                                & \text{ else. }
\end{cases}
\]
and 
\[  
i_S^!\mcH^k( IC_t[d])= \mcH^{k+d}(\mathbb{D}_S i_S^* IC_{t^*})= i_S^*\mcH^{-k-d-2d_S}(IC_{t^*})
\]
implies 
\[  
i_S^!\mcH^k( IC_t[d])= 
\begin{cases}
%0,                               & \text{ if }  d_S  > d_{S^\prime} \\
%                                 & \text{ or }  d_S =  d_{S^\prime}  , i\neq d_S \\
%                                 & \text{ or }  d_S  < d_{S^\prime},   i\in (-\infty ,-d_{S^\prime}-1]\cup [-d_S, \infty)    \\
\mcL^* ,                       & \text{ if }  d_S =  d_{S^\prime} ,  k+d =-d_S \\
\Ind_{S^\prime}^S (\mcL^*)_{-k-d-2d_S} & \text{ if }  d_S  < d_{S^\prime},   -k-d-2d_S\in [-d_{S^\prime}, -d_S-1]  \\
0                                & \text{ else. }
\end{cases}
\]
where $d_S= \dim_{\C}S$. This implies 
\[ 
\begin{aligned}
H^{[k]}(\pi^{-1}(s)) &= \bigoplus_t \bigoplus_{d\in \Z} L_{t,d} \otimes_{\C} i_S^* \mcH^{k+d} (IC_t) \\
                                           &= \bigoplus_{t=(S,\mcL)}L_{t, -d_s-k} \otimes_{\C} \mcL  \oplus  
                                              \underbrace{\bigoplus_{t=(S^\prime , \mcL), d_S < d_{S^\prime}}  \; 
                                              \bigoplus_{r=-d_{S^\prime}}^{-d_S-1} L_{t,r-k}\otimes_{\C} \Ind_{S^\prime}^S (\mcL)_r}_{=:H^{[k]}(\pi^{-1}(s))_{>S}}
\end{aligned}
\]
and
\[ 
\begin{aligned}
H_{[k]}(\pi^{-1}(s)) &= \bigoplus_t \bigoplus_{d\in \Z} L_{t,d} \otimes_{\C} i_S^! \mcH^{i+d} (IC_{t^*}) \\
                                           &= \bigoplus_{t=(S,\mcL)}L_{t, -d_s-k} \otimes_{\C} \mcL^*  \oplus  
                                              \underbrace{\bigoplus_{t=(S^\prime , \mcL), d_S < d_{S^\prime}}  \; 
                                              \bigoplus_{r=-d_{S^\prime}}^{-d_S-1} L_{t,-r-2d_S-k}\otimes_{\C} \Ind_{S^\prime}^S (\mcL^*)_r}_{=:H_{[k]}(\pi^{-1}(s))_{>S}}
\end{aligned}
\]
as $\pi_1(S,s)$-representations. 
We call the direct summands isomorphic to $\Ind_{S^\prime}^S (\mcL)_r, r\in [-d_{S^\prime}, -d_S-1] $ the \emph{unwanted summands}. 
Now we can explain how you can recover from the $\pi_1(S,s)$-representations $H^{[k]}(\pi^{-1}(s)), k\in \Z$ the 
data for the decomposition theorem (i.e. the local sytems and the graded multiplicity spaces). If $d_S$ is the maximal one, it holds 
\[H^{[*]}(\pi^{-1}(s)) =\bigoplus_{k\in \Z} \bigoplus_{t=(S,\mcL)}L_{t, -d_s-k} \otimes_{\C} \mcL  \]
and we can recover the graded multiplicity spaces $L_t$ with $t=(S,?)$ for the dense stratum ocurring in the decomposition theorem. 
For arbitrary $S$ we consider 
\[ 
H^{[*]}(\pi^{-1}(s)) / H^{[*]}(\pi^{-1}(s))_{>S}  \cong \bigoplus_{k\in \Z} \bigoplus_{t=(S,\mcL)}L_{t, -d_s-k} \otimes_{\C} \mcL 
\]
and by induction hypothesis we know the $\pi_1(S,s)$-representation $H^{[*]}(\pi^{-1}(s))_{>S}$, therefore we can recover the $L_t$ with $t=(S,?)$ 
from the above representation. 

Now assume that $\pi$ is semi-small. % (i.e. $\dim Z_{i,i}=e_i, \; i\in I$). 
Then, we know that $L_{t,d}=0$ for all $t=(S,\mcL) $ whenever $d\neq 0$. We can also restrict our attention on a direct summand 
$(\pi_i)_*\underline{\C}[e_i]$ for one $i\in I$ and find the decomposition into simple perverse sheaves. 
That means we only need $H^{e_i-d_S}(\pi_i^{-1}(s)$ to recover the data for the decomposition theorem. 
It also holds $2\dim \pi_i^{-1}(s) \leq e_i-d_S, i\in I$ and since $H^{e_i-d_S}(\pi_i^{-1}(s))=0$ whenever $2\dim \pi_i^{-1}(s)< e_i-d_S$, 
we only need to consider the strata $S$ with $2\dim \pi_i^{-1}(s) = e_i-d_S$, then 
\[ H^{e_i-d_S}(\pi_i^{-1}(s))= H^{top}(\pi_i^{-1}(s)) \neq 0\]
and we call $S$ a \textbf{relevant stratum} for $i$  $(\in I)$. We call a stratum relevant if it is relevant for at most one $i\in I$. \\
Analogously, one can replace $H^{[k]}(\pi^{-1}(s))$ by $H_{[-k]}(\pi^{-1}(s))$ and stalk by costalk. \\

Let $x\in V$ be arbitrary. By a previous section we know that $H_{[*]}(\pi^{-1}(x))$ and $H^{[*]}(\pi^{-1}(x))$ are left (and right) 
graded $H_{[*]}(Z)$-modules.
The following lemma explains their special role. Unfortunately, the following statement is only known if all 
strata $S$ contain a $G$-orbit $s\in \mcO \subset S$ such that $\pi_1(\mcO ,s) =\pi_1 (S,s)$. For local systems on the strata this is by monodromy 
the same as the assumption that all strata are $G$-orbits. Let $C$ be a finite group, we write $\Simp(C)$ for the set of isomorphism classes of simple 
$\C C$-modules and denote by $\bf{1}\in \Simp(C)$ the trivial representation\footnote{In the literature this is called $Irr(C)$, we use the word irreducible 
only for a property of topological spaces}.  

\begin{lemma} (\cite{CG}, Lemma 8.4.11, p.436, Lemma 3.5.3, p.170) 
Assume that the image of the Springer map contains only finitely many $G$-orbits. 
\begin{itemize}
\item[(a)] Let $\mcO=Gx\subset V$ be a $G$-orbit. There is an equivalence of categories between 
\[\{\text{$G$-equivariant local systems on }\mcO\} \leftrightarrow C(x)\!-\!\Mod\]  
where $C(x)=\Stab_G(x)/(\Stab_G(x))^o$ is the \textbf{component group} of the stabilizer of $x$. In particular, via monodromy also 
the $\pi_1(\mcO, s)$-representations which correspond to $G$-equivariant local systems on $\mcO$ are equivalent to $C(x)\!-\!\Mod$. 
\item[(b)] The $C(x)$-operation and the $H_{[*]}(Z)$-operation on $H_{[*]}(\pi^{-1}(x))$ (and on $ H^{[*]}(\pi^{-1}(x))$) commute. 
\end{itemize}
\end{lemma}

The semi-simplicity of $C(x)\!-\!\Mod$ implies that
\[
H_{[*]}(\pi^{-1}(x)) = \bigoplus_{k\in \Z}\bigoplus_{\chi \in \Simp(C(x))} (H_{[k]}(\pi^{-1}(x)))_{\chi} \otimes_{\C} \chi     
\] 
where $\Simp(C(x))$ is the set of isomorphism classes of simple $C(x)$-modules and for any $C(x)$-module $M$ we call 
$M_{\chi }:=\Hom_{C(x)\!-\!\Mod} (\chi ,M)$ an \textbf{isotypic component}. 
Since the two operation commute it holds $(H_{[*]}(\pi^{-1}(s))_{\chi }$ naturally has the structure of a graded $H_{[*]}(Z)$-module. But we will from now just see it as a module over $H_*(Z)$.  
As $H_{*}(Z)-C(x)$-bimodule decomposition we can write the previous decomposition as 
\[
H_{[*]}(\pi^{-1}(s)) = \bigoplus_{\chi \in \Simp(C(x))} H_*(\pi^{-1})_{\chi} \boxtimes \chi      
\] 
where $H_*(\pi^{-1})_{\chi} \boxtimes \chi$ is the obvious bimodule $H_*(\pi^{-1}(x))_{\chi}\otimes \chi$. As an immediate consequence of this we get, 
if $Gx$ is a dense orbit in the image of the Springer map, then 
\[
L_{t,-*}(-d_{Gx}) = H_{[*]}(\pi^{-1}(x))_{\chi}, \text{ for }t=(x, \chi), \chi \in \Simp(C(x)),
\]
in particular, $H_{[*]}(\pi^{-1}(x))$ is a semisimple $H_{*}(Z)$-module (graded and not graded), even a semisimple $H_{*}(Z)-C(x)$-bimodule. 
For more general orbits, we do not know if it is semi-simple. In the case of a semi-small Springer map we have the following result. 

\begin{satz} \label{SpringerCorres} Assume the image of the Springer map $\pi$ has only finitely may orbits 
and $\pi$ is semi-small. 
%\begin{itemize}
%\item[(i)] All (left and right) Springer fiber modules $H_*(\pi^{-1}(x)), H^*(\pi^{-1}(x), x\in V$ are semi-simple $H_*(Z)$-modules. Every simple 
%$H_*(Z)$ module occur as a dirct summand of a Springer fibre module. 
There is a bijection between the following sets 
\begin{itemize}
\item[(1)] 
$\{ (x,\chi ) \mid \mcO=Gx , \chi \in \Simp (C(x)),  H_{[d_{\mcO}]}(\pi^{-1}(x))_{\chi}\neq 0\}$ where the $x$ in $V$ are 
in a finite set of points representing the $G$-orbits in the image of the Springer map. 
\item[(2)] $\Simp (H_{<0>}(Z)\! -\!\Mod):=$ simple $H_{<0>}(Z)$-modules up to isomorphism 
\item[(3)] $\Simp (H_{<*>}^A(Z)\!-\!\Mod^{\Z})$ := simple graded $H_{<*>}^A(Z)$-modules up to isomorphism and shift for any $A\in \{pt, T, G\}$. 
\end{itemize}
Between (1) and (2), it is given by $(x,\chi)\mapsto H_{[d_{\mcO}]}(\pi^{-1}(x))_{\chi}$. 
We call this bijection the \textbf{Springer correspondence}. \\
For a relevant orbit $\mcO$ (for at least one $i\in I$) it holds 
\[H_{[d_{\mcO}]}(\pi^{-1}(x))_{\bf{1}}=\bigoplus_{i\colon 2\dim \pi_i^{-1}(s) = e_i-d_{\mcO}} H_{top}(\pi_i^{-1}(x))^{C(x)}\neq 0\] 
and $C(x)$ operates on the top-dimensional irreducible components of $\pi_i^{-1}(x)$ by permutation. This implies we get an injection 
\[
\begin{aligned}
\{ \text{relevant $G$-orbits in }Im(\pi)\}  &\inj \Simp (H_{<0>}(Z)\!-\!\Mod)  \\  
 \mcO=Gx &\mapsto H_{[d_{\mcO}]}(\pi^{-1}(x))^{C(x)} 
\end{aligned}
 \]
\end{satz}

\paragraph{sketch of proof:} For $k=d_{\mcO}$ look at the decomposition for $H_{[k]}(\pi^{-1}(x)$ and use that $L_{t,d}=0$ whenever $d\neq 0$ to see that 
the unwanted summands vanish. Then show that the decomposition coincides with the second decomposition (with respect to the irreducible characters of $C(x)$) of  $H_{[k]}(\pi^{-1}(x))$ which gives the identification 
of the $L_t$ with the $H_{[d_{\mcO}]}(\pi^{-1}(x))_{\chi}$. \hfill $\Box$  
%Set $Z_{\mcO }:= p^{-1}(\mcO ) =\pi^{-1}(\mcO )\times_{\mcO} \pi^{-1}(\mcO )\cong G\times^{Stab (x)} [\pi^{-1}(x)\times \pi^{-1}(x)]$, there is an isomorphism 
%\[ H_{[d_{\mcO}]}(Z_{\mcO}) \to H_{[d_{\mcO}]}(\pi^{-1}(x)\times \pi^{-1}(x))^{C(x)} \]
%of $H_{<0>}(Z)$-bimodules. Choose a total order $\leq$ refining the partial order in the orbits given given by inclusion into the closure. Set $Z_{\leq \mcO}=\bigcup_{\mcO^\prime \leq \mcO}$, we get a filtration of 2-sided ideals $H_{[--]}(Z_{\leq \mcO})$ in $H_{<0>}(Z)$. Since the algebra is semi-simple we have get an isomorphism of $H_{<0>}(Z)$-bimodules 
%\[
%\begin{aligned}
%H_{<0>}(Z)\cong gr_{\mcO} H_{<0>}(Z) &=\bigoplus_{\mcO}H_{[--]}(Z_{\mcO})=\bigoplus_{\mcO=Gx}\Hom_{C(x)-mod}(H_{[d_{\mcO}]}(\pi^{-1}(x)),H_{[d_{\mcO}]}(\pi^{-1}(x))) \\
%&= \bigoplus_{\mcO}\bigoplus_{\chi, \psi \in Simp(C(x))} \Hom_{C(x)} (\chi, \psi ) \otimes \Hom_{\C}( H_{[d_{\mcO}]}(\pi^{-1}(x))_{\chi},  H_{[d_{\mcO}]}(\pi^{-1}(x))_{\psi}) \\
%&= \bigoplus_{\mcO=Gx, \chi \in Simp(C(x))} \End_{\C} ( H_{[d_{\mcO}]}(\pi^{-1}(x))_{\chi})
%\end{aligned}
%\]
%where for the last equality one uses $\Hom_{C(x)}(\chi , \psi)=\delta_{\chi , \psi} \C$. This implies the claim. \hfill $\Box$

It is an open question to understand Springer fibre modules more generally. Also, Springer correspondence hints at a hidden equivalence of categories. 
This functorial point of view we investigate in the next subsection. 

%%%%%%%%%%%%%%%%%%%%%%%%%%%%%%%%%%%%%%%%%%%%%%%%%%%%%%%%%%%%%%%%%%%%%%%%%%%%%%%%%%%%%%%%%%%%%%%%%%%%%%%%%%%%%%%%%%%%%%%%%%
%%%%%%%%%%%%%%%%%%%%%%%%%%%%%%%%%%%%%%%%%%%%%%%%%%%%%%%%%%%%%%%%%%%%%%%%%%%%%%%%%%%%%%%%%%%%%%%%%%%%%%%%%%%%%%%%%%%%%%%%%%

\subsection*{The Springer functor}
We consider $H_{[*]}^A(Z)$ again with the grading from the thm \ref{Extalgebra}.
Let $proj^{\Z} H_*^A(Z)$ be the category of finitely generated projective $\Z$-graded left $H_{[*]}^A(Z)$-modules, 
morphisms are the module homomorphisms which 
are homogeneous of a degree $0$. 
Let $\mcP^A\subset D_A^b(X)$ be the full subcategory closed under direct sums and shifts generated by $IC_t^A$, $t=(S, \mcL)$ be the tuple of 
strata with simple local system on it which occur in the decomposition theorem.   

The following lemma is in a special case due to Catharina Stroppel and Ben Webster, see \cite{SW}.   
\begin{lemma}\label{Springer-functor}
The functor 
\[
\begin{aligned}
proj^{\Z} H_{[*]}^A(Z) &\to \mcP_A \\
             M&\mapsto \bigoplus_{i\in I} (\pi_i)_*\underline{\C}_{E_i}[e_i] \otimes_{H_{[*]}^A(Z)} M
\end{aligned}
\]
is an equivalence of semisimple categories mapping $P_t^A \mapsto IC_t^A$. We call this the \textbf{Springer functor}
\footnote{This name is due to Dustin Clausen in his thesis.}
\end{lemma}

\paragraph{proof:} 
By thm \ref{Extalgebra} we know $H_{[*]}^A(Z) = 
\Ext_{D_A^b(V)}^*(\bigoplus_{i\in I} (\pi_i)_*\underline{\C}_{E_i}[e_i],\bigoplus_{i\in I} (\pi_i)_*\underline{\C}_{E_i}[e_i])$ is 
an isomorphism of graded algebras. This makes the functor well-defined. 
The direct sum decomposition of $\bigoplus_{i\in I} (\pi_i)_*\underline{\C}_{E_i}[e_i] $ by the decomposition theorem in $\mcP_A$ corresponds to idempotent 
elements in $H_{[0]}^A(Z)$, 
which correspond (up to isomorphism and shift) to the indecomposable projective graded modules, let for example $P_t = H_{[*]}^A(Z) e_t$. 
Shifts of graded modules are mapped to shifts in $\mcP_A$, 
therefore the functor is essentially surjective. It is fully faithful because of the mentioned isomorphism 
\[ 
\Hom_{proj^{\Z} H_{[*]}^A(Z)}(P_t, P_s (n)) = e_s H_{[n]}^A(Z) e_t = \Hom_{D_A^b(V)}( IC_t , IC_s[n]) 
\]
\hfill $\Box$

Let $\mcP^A(V) \subset D_A^b(V)$ be the category of $A$-equivariant perverse sheaves on $V$. Assume for a moment that the map $\pi$ is semi-small. 
Then, we know that $\bigoplus_{i\in I} (\pi_i)_*\underline{\C}_{E_i}[e_i] $ is an object of $\mcP^A(V)$. In this situation the two gradings of the 
Steinberg coincide. The top-dimensional Borel-Moore homology $H_{top}(Z_{i,i})$ 
coincides with the degree zero subalgebra $H_{[0]}(Z_{i,i})$. We want the Springer 
functor to go to a category of perverse sheaves, i.e. we do not want to allow shifts of the grading for modules. 
Therefore, we pass to 
\[ H_{[0]}(Z)=H_{<*>}(Z)/(H_{<*>}(Z))_{>0} = H_{<*>}^A(Z)/(H^A_{<*>}(Z))_{>0}, \quad A\in \{ pt, T, G\}\]
and replace projective graded modules over $H_{[*]}^A(Z)$ by the additive category of simple modules over $H_{[0]}(Z)$, this equals the 
category $H_{[0]}(Z)\!-\!\Mod$ of finite dimensional (ungraded) modules over $H_{[0]}(Z)$ because the algebra is semi-simple. \\
In particular, it holds \\
$H_{[0]}(Z) =\Ext^0_{D^b_A(V)}( \bigoplus_{i\in I} (\pi_i)_*\underline{\C}_{E_i}[d_i] , \bigoplus_{i\in I} (\pi_i)_*\underline{\C}_{E_i}[d_i] )
=\End_{\mcP^A (V)}( \bigoplus_{i\in I} (\pi_i)_*\underline{\C}_{E_i}[d_i] ), $  $A\in \{pt, T, G\}$. 

The following lemma is for classical Springer Theory due to Dustin Clausen, cp. Thm 1.2 in \cite{Cl}. 
\begin{lemma} If the Springer map $\pi$ is semi-small, 
we have the following version of the Springer functor 
\[
\begin{aligned}
\mcS\colon H_{[0]}(Z)\! -\! \Mod & \to \mcP^G(V) \\
M & \mapsto 
\bigoplus_{i\in I} (\pi_i)_*\underline{\C}_{E_i}[e_i] \otimes_{H_{[0]}(Z)} M
\end{aligned}
\]
It holds that $\mcS$ is an exact functor (between abelian categories) and it is fully faithful. 
If $e_i+e_j$ is even for all $i,j\in I$ the $\mcS$ identifies $H_{[0]}(Z)\!-\!\Mod$ with a semi-simple Serre subcategory of 
$\mcP^G(V)$ (i.e. it is an exact subcategory which is also extension closed). 
Furthermore it is invariant under Verdier duality on $\mcP^G(V)$. 
\end{lemma}

\begin{rem}
Assume that the Springer map is semi-small, 
the image of the Springer map contains only finitely many $G$-orbits and each $G$-orbit is relavant and simply connected, then the Springer functor 
from above is an equivalence of categories. (The only known example for this is the classical Springer map for $G=\Gl_n$, see later.)  
\end{rem}

\paragraph{proof:}
A similar proof as in the lemma above shows that the Springer functor induces an equivalence on the 
full subcategory of $\mcP^G(V)$ generated by finite direct sums of direct summands of $ \bigoplus_{i\in I} (\pi_i)_*\underline{\C}_{E_i}[e_i]$. 
This is a semi-simple category. 
Assume that $e_i+e_j$ is even for all $i,j\in I$, we have to see that it is extension closed.  
By composition with the forgetful functor we get a functor 
\[ H_{[0]}(Z)\!-\!\Mod  \xrightarrow{\mcS} \mcP^G(V) \xrightarrow{F} \mcP^{pt}(V) =: \mcP (V), \]
by \cite{Cl} the forgetful functor $F$ is fully faithful. Now, by \cite{Ar3}, 4.2.10 
the category $\mcP (V)$ of $D^b(V)$ is closed under extensions and admissible because it is the heart of a t-structure. By the Riemann Hilbert correspondence 
there exists an abelian category $\mcA$ ($=$ regular holonomic $D$-modules on $V$) and an equivalence of triangulated categories ($=$ the de Rham functor)
\[ DR_V\colon D^b(\mcA) \to D^b(V) \]
such that the standard $t$-structure on $D^b(\mcA)$ is mapped to the perverse $t$-structure and it restricts to an equivalence of categories 
$\mcA \to \mcP (V)$. This implies that for $X\cong DR_V(X^\prime), Y\cong DR_V(Y^\prime) $ in $\mcP (V)$ and $n\in \N_0$  
\[
\Ext^n_{\mcP(V)} (X,Y) = \Ext^n_{\mcA}(X^\prime, Y^\prime)  =\Hom_{D^b(\mcA)}(X^\prime , Y^\prime [n]) = \Hom_{D^b(V)} (X, Y[n]) 
\]
where the first and the third equality follows from the de Rham functor and the second equality holds because it is the standard $t$-structure, 
cp. for example \cite{GM}, p.286. 

Now, since we know 
\[
\Hom_{D^b(V)}( \bigoplus_{i\in I} (\pi_i)_*\underline{\C}_{E_i}[e_i], (\bigoplus_{i\in I} (\pi_i)_*\underline{\C}_{E_i}[e_i])[1]) 
=H_{<1>}(Z) =\bigoplus_{i,j\in I} H_{e_i+e_j-1}(Z) =0
\]
because $H_{odd}(Z)=0$ by lemma \ref{oddVan} and the assumption that $e_i+e_j$ is even for every $i,j\in I$. We obtain that 
\[\Ext^1_{\mcP (V)} ( \bigoplus_{i\in I} (\pi_i)_*\underline{\C}_{E_i}[e_i], \bigoplus_{i\in I} (\pi_i)_*\underline{\C}_{E_i}[e_i]) =0, \] 
i.e. the semi-simple category generated by the direct image of the Springer map is extension closed. 
\hfill $\Box $

\subsection*{ What is Springer Theory ?}
%We can summarize this (roughly) to the following:
One possible definition 
\[
\boxed{
\begin{aligned}
\text{Springer theory (for }& (G,P_i,V,F_i)_{i\in I} \text{ and a choice of }H \text{) is to understand }\\
\text{ the Steinberg algebra }& \text{ together with its graded modules.}  
\end{aligned}
}
\]

But I think today it is sensible to say Springer theory is the study of all categories and algebras (and modules 
over it) which have a construction originating in some Springer Theory data $(G,P_i, V, F_i)_{i\in I}$. Then, this includes 
  
\begin{itemize}
\item[(1)] Monoidal catgories coming from multiplicative families of Steinberg algebras and their Grothendieck ring. In particular, this includes 
Lusztig's categories of perverse sheaves (see \cite{L} and the example quiver-graded Springer Theory later). 
\item[(2)] Noncommutative resolutions\footnote{\emph{here}: This means just a tilting vector bundle on $E$, 
because this gives $t$-structures in the category of coherent sheaves on $E$} corresponding to the 
Springer map. In particular, this includes Bezrukavnikov's noncommutative counterparts of the Springer map in 
\cite{Be} and Buchweitz, Leuschke and van den Bergh's articles \cite{BLvdB} and \cite{BLvdB2}. 
\item[(3)] Categories of flags of ($KQ$-)submodules for given quivers because their isomorphism classes parametrize orbits of (quiver-graded) 
Springer fibres. This includes for example Ringel's and Zhang's work on submodule categories and preprojective algebras \cite{Ri4}. 
Also certain $\Delta$-filtered modules studied in  \cite{BHRR}, \cite{BH}.  
\end{itemize}

An (of course) incomplete overview can be found in the flowchart at the end of this article. 

We would also like you to observe that in the two examples we explore connections 
between objects roughly related to the following triangle 

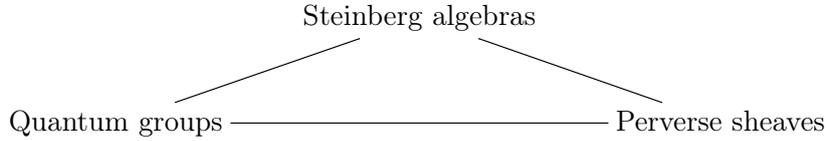
\[
\xymatrix{
&\text{Steinberg algebras} \ar@{-}[dr] \ar@{-}[ld] & \\
\text{Quantum groups} && \text{Perverse sheaves} \ar@{-}[ll]
}
\]

\subsection*{Classical Springer Theory }

This is the case of the following initial data
\[
\left[
\begin{aligned}
(*)&\; G\text{ an arbitrary reductive group,}\\ 
(*)&\; P=B\text{ a Borel subgroup of }G\text{, denote its Levi decomposition by }B=TU\\
   &\text{ with }T\text{ maximal torus, }U\text{ unipotent.} \\
(*)&\; V=\mathfrak{g}\text{ the adjoint representation,}\\
(*)&\; F=\mathfrak{n}:= \Lie (U).  
\end{aligned}
\right.
\] 
We set $\mathcal{N}:= G\mathfrak{n}$, i.e. the image of the Springer map, and call it the nilpotent cone. We consider the Springer map 
as $\pi \colon E=G\times^B \mathfrak{n} \to \mathcal{N}$. Explicitly, we can write the Springer triple as 
\[
\xymatrix{  &E= \{ (n,gB) \in \mcN \times G/B \mid n\in {}^g\mathfrak{b}:= \Lie(gBg^{-1})\} \ar[dl]^{\pi=pr_1} \ar[dr]_{\mu =pr_2} & \\
\mcN && G/B 
}
\]
For $G=\Gl_n$ we identify $\Gl_n/B$ with the variety $Fl_n$ if complete flags in $\C^n$ and 
\[E=\{ (A, U^\bullet ) \in \End_{\C}(\C^n) \times Fl_n \mid 
A^n=0, A(U^k)\subset U^k, \; 1\leq k\leq n \}.\] 
It turns out, $\pi$ can be identified with the moment map of $G$, in particular,  
$E\cong T^*(G/B)$ is the cotangent bundle over $G/B$ and $\pi$ is a resolution of singularities for $\mathcal{N}$. But most importantly, this makes the Springer map a symplectic resolution of singularities and one can use 
symplectic geometry (see for example \cite{CG}). \\
The Steinberg variety is given by 
\[
\xymatrix{  &Z= \{ (n,gB, hB) \in \mcN \times G/B \times G/B \mid n\in {}^g \mathfrak{b} \cap{}^h \mathfrak{b}  \} \ar[dl]^{p=pr_1} \ar[dr]_{m=pr_{2,3}} & \\
\mcN && G/B \times G/B 
}
\]
for $G=\Gl_n$ we can write it as 
\[Z=\{ (A, U^\bullet ,V^\bullet ) \in \End_{\C}(\C^n) \times Fl_n \times Fl_n \mid 
A^n=0, A(U^k)\subset U^k , A(V^k)\subset V^k, \; 1\leq k\leq n  \}.\] 
Recall, that we had the stratification by relative position $Z^w:= m^{-1}( G\cdot (eB, wB)), w\in W$ where $W$ is the Weyl group of $G$ with respect to a maximal torus $T\subset B$. Since $Z^w\to G\cdot (eB, wB)$ is a vector bundle, we can easily calculate its dimension 
\[ 
\begin{aligned}
\dim Z^w &= \dim G\cdot (eB, wB) + \dim \mathfrak{n}\cap {}^w \mathfrak{n} \\
         &= \dim G - \dim B\cap {}^wB +\dim \mathfrak{n}\cap {}^w \mathfrak{n}  = \dim G -\dim T \\
         &= \dim E.
\end{aligned}
\]
We conclude that $Z$ is equidimensional of dimension $e:= \dim E$, in particular the Springer map is semi-small. Also we see that 
the irreducible components of $Z$ are given by $\overline{Z^w}, w\in W$, that implies that the top-dimensional Borel-Moore homology group $H_{top}(Z)$ has a $\C$-vector space  basis given by the cycles $[\overline{Z^w}]$. 
In the semi-small case we know $H_{[0]}(Z)=H_{<0>}(Z)=H_{top}(Z)$ is a sub- and quotient algebra of $H_*(Z)$. 

\begin{exa} $G=\Sl_2$, $B=\{ \left(\begin{smallmatrix} a & b \\ 0 & a^{-1} \end{smallmatrix}\right) \mid b \in \C, a\in \C\setminus \{ 0\}\}$. Then $\mcN= \{ (x,y,z) \in \C^3 \mid x^2+zy =0\}$ and \[E = \{ (A, L) \in \mcN\times \mathbb{P}^1\mid L \subset \Ker A \} = \{ (\left(\begin{smallmatrix} x & y \\ z & -x \end{smallmatrix}
\right), [a\colon b]) \in M_2 (\C )\times \mathbb{P}^1 \mid \; x^2+zy=0, xa +yb=0, za-yb=0\},\]
the Springer map can be seen as the following picture 
 
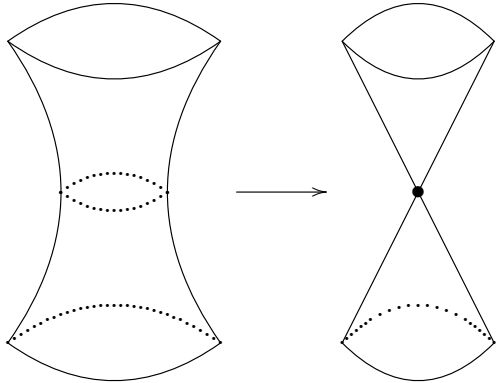
\[
 \xy
    (0,0)*++{\xy
  (0,0)*{}="A";
    (28,0)*{}="B";
    (0,40)*{}="C";
    (28,40)*{}="D";
    (7,20)*{}="E";
    (21,20)*{}="F";
 "A"; "C" **\crv{(14,20)};
 "C"; "D" **\crv{(14,30)};
 "C"; "D" **\crv{(14,50)};
"B"; "D" **\crv{(14,20)};
"A"; "B" **\crv{(14,-10)};
 "A"; "B" **\crv{~*=<\jot>{.}(14,10)};
 "E"; "F" **\crv{~*=<\jot>{.}(14,15)};
 "E"; "F" **\crv{~*=<\jot>{.}(14,25)};
    %"A"; "D" **\dir{-};
    %"B"; "C" **\dir{-};    
\endxy }="x";
(40,0)*++{\xy 
 (0,0)*{}="A";
    (20,0)*{}="B";
    (0,40)*{}="C";
    (20,40)*{}="D";
 "C"; "D" **\crv{(10,30)};
 "C"; "D" **\crv{(10,50)};
    "A"; "D" **\dir{-};
    "B"; "C" **\dir{-};
 "A"; "B" **\crv{(10,-10)};
"A"; "B" **\crv{~*=<\jot>{.}(10,10)};
(10,20)*{\bullet};
\endxy }="y";
    {\ar "x";"y"};
    \endxy
    \]

This is well-known to be the crepant resolution of the $A_2$-singularity from the MacKay correspondence. 
In general, if $G$ is semi-simple of type $ADE$, then there exists a slice of the nilpotent 
cone such that the restricted map is the crepant resolution of the corresponing type singularity, see \cite{Sl} for more details. 
\end{exa}

\begin{satz} (roughly  Springer \cite{Spr1}) 
There is an isomorphism of $\C$-algebras 
\[
\begin{aligned} 
H_{top} (Z) &\cong \C[W] \\
[\overline{Z^s}] &\mapsto s-1 
\end{aligned}
\]
\end{satz}
%and $H^{top}(Z)$ is a semi-simple algebra, it is isomorphic to $H_{*}(Z)/(H_{*}(Z))_{>0}$. The semi-smallness of the Springer map says that 
%$\pi_* \underline{C}[d]$, $d=\dim_{\C} E$ is a direct sum of simple perverse sheaves, it follows that we write down a version of the Springer functor 
%due to D. Clausen (in \cite{Cl}) which maps to the category of perverse sheaves $\mcP (\mathcal{N})$ (we do not need the shifts). Since we work now with the 
%simple top of the algebra we have to replace the category of projective modules over the Steinberg algebra by the semi-simple category generated by their simple tops, i.e. by $\C W-mod$. Then, the Springer functor  
The Springer functor (due to Clausen, \cite{Cl}) takes the form 
\[
\begin{aligned}
\C W\! -\!\Mod & \to \mcP^G (\mathcal{N})\\
M & \mapsto \pi_*\underline{\C}[e] \otimes_{\C W} M
\end{aligned}
\]
and identifies $\C W\! -\!\Mod$ with a semi-simple Serre subcategory of $\mcP^G (\mcN )$. This implies an injection on simple 
objects which are in $\mcP^G(\mcN )$ the intersection cohomology complexes associated to $(\mcO , \mcL)$ with $\mcL$ a 
simple $G$-equivariant local system on an $G$-orbit $\mcO\subset \mcN$. As a consequence we get the bijection called \emph{Springer correspondence} from 
thm \ref{SpringerCorres}  
\[ 
\begin{aligned}
\Simp (W) &\Leftrightarrow  \{ t= (\mcO, \mcL) \mid  \text{ certain (=occurring in the decomp. thm) }\}\\
         &= \{ (x, \chi) \mid x \in \mcN \text{ rep of $G$-orbits }, \chi \in \Simp (C(x)), (H_{top} (\pi^{-1}(x)))_{\chi}\neq 0 \}
\end{aligned}         
\]
where $\Simp(W)$ is the set of isomorphism classes of simple objects in $\C W\! -\!\Mod $. 
The inverse of the map is given by $(x,\chi)\mapsto (H_{top} (\pi^{-1}(x)))_{\chi}$. 
For this Springer map all orbits in $\mcN $ are relevant, i.e. we also have an injection 
\[  
\begin{aligned}
\{ \text{$G$-orbits in }\mcN \} & \to \Simp (W) \\
        Gx &\mapsto H_{top}(\pi^{-1}(x))^{C(x)}
\end{aligned}        
\]

%The theorem \ref{SpringerCor} is the main result of classical 
%Springer theory, we do not repreat the statement here. 

%in the words of proposition .. we get the following  
%\begin{satz} 
%There is a bijection 
%\[
%\begin{aligned}
%\{ (x, \chi) \mid 
%x\in \mcN \text{ rep of $G$-orbits}, \chi \in Simp (C(x)), H_{top}(\pi^{-1}(x))_{\chi}\neq 0 \} &\to Simp (W)  \\
%(x, \chi ) & \mapsto H_{top}(\pi^{-1}(x))_{\chi}
%\end{aligned}
%\]
%which is called \textbf{Springer correspondence}. It always holds $H_{top}(\pi^{-1}(x))_{1}=H_{top}(\pi^{-1}(x))^{C(x)}\neq 0$ because $C(x)$ operates on the top-dimensional irreducible components of $\pi^{-1}(x)$ by permutation. This implies we get an injection 
%\[
%\begin{aligned}
%\{ \text{$G$-orbits in }\mcN\}  &\inj Simp (W)  \\  
% \mcO=Gx &\mapsto H_{top}(\pi^{-1}(x))^{C(x)} 
%\end{aligned}
% \]
%\end{satz}

\begin{rem} \label{others} We remark that there are several alternative constructions of the group operation of $W$ on the  
Borel-Moore homology/ singular cohomology of the Springer fibres. In \cite{Ar3}, section 5.5 you find an understandable treatment of Lusztig's approach to this 
operation using intermediate extensions for perverse sheaves and Arabia provides a list of other authors and approaches to this (first Springer \cite{Spr1},\cite{Spr2}, then Kazhdan-Lusztig \cite{KL2}, Slodowy \cite{Sl2}, Lusztig \cite{L2}, Rossmann \cite{Ros}) and these operations differ between each other by at most by multiplication with a sign character (see \cite{Ho}). \\
Also, Springer proves with taking (co)homology of Springer fibres with rational coefficients that the simple $W$-representations are 
all even defined over $\Q$, a result which our approach does not give because the 
simple $C(x)$-modules are not necessarily all defined over $\Q$ (cp. \cite{CG}, section 3.5, p.170). In Carter's book \cite{Ca}, p. 388, you find for simple adjoint groups the component groups $C(x), x\in \mcN$ are one of the following list $ (\Z/2\Z)^r, \;\; S_3, S_4, S_5, \quad r\in \N_0 $ 
as a consequence he gets that the simple modules over the group ring are already defined over $\Q$. 
\end{rem}

%
%\begin{rem}
%Lusztig's generalized Springer theory. Brylinsky?? realization of Springer correspondence via Fourier transform of perverse sheaves. 
%Clausen: Springer functor which induces the bijection between simples, Grinberg 98: Generalized Springer correspondence. 
%\end{rem}
In the introduction of the book \cite{BBM} you find for a semisimple group $G$ a triangle 
\[
\xymatrix{
& \text{simple $\C W$-modules} & \\
\text{primitive ideals in }U(\mathfrak{g})\ar[ur]\ar[rr] && \text{$G$-orbits in the nilpotent cone }\ar[ul] 
}
\]

They explain it as follows (i.e. this is a summary of a their summary). 
\begin{itemize}
\item[*] There is an injection of $G$-orbits in $\mcN$ into simple $\C W$-modules by the Springer correspondence. 
\item[*] A primitive ideal in $U(\mathfrak{g})$ is a kernel of some simple $U(\mathfrak{g})$-representation. The classification of 
primitive ideals is archieved as a result of the proof of the Kazhdan-Lusztig conjectures (see Beilinson-Bernstein \cite{BB}, Brylinski-Kashiwara \cite{BK}). Any ideal in 
$U(\mathfrak{g})$ has an associated subvariety of $\mathfrak{g}$. The associated variety of a primitive ideal is the closure of an orbit in $\mcN$, this was first conjectured by Borho and Jantzen. 
\item[*] Joseph associated to a primitive ideal a $W$-harmonic polynomial in $\C[\mathfrak{t}]$ (=Goldie rank polynomial) which is a basis element of one of the 
simple $\C W$-modules. 
\end{itemize}

We also have to mention the following important results which use $K$-theory instead of Borel-Moore homology. 
\paragraph{Parametrizing simple modules over Hecke algebras.}
This field goes back to the work of Kazhdan and Lusztig on the proof of the Deligne-Langlands conjecture for Hecke algebras, see \cite{KL87}. 
They realize irreducible representations of Iwahori Hecke algebra as Grothendieck groups of equivariant (with respect to certain groups) coherent 
sheaves on the  Springer fibres. This is now known as Deligne-Langlands correspondence and we call similar results which come later for different Hecke algebras still DL-correspondence.\\ 
\noindent 
Let $G$ be an algebraic group and $X$ a $G$-variety, let $K_0^G(X):=K_0(coh^G(X))$ be the Grothendieck group of the category of $G$-equivariant coherent sheaves on $X$. The group $\C^*$ operates on the (classical) Steinberg variety $Z$ via $(n, gB, hB)\cdot t:= (t^{-1}n, gB,hB)$, 
the convolution product construction gives a ring structure on 
$K_0^{G\times \C^*}(Z)$. \\
Recall, for a reductive group we fix 
a maximal torus and a Borel subgroup 
$T\subset B\subset  G$ and call $(W,S)$ the associated Weyl group with set of simple roots. We write $X(T)=\Hom (T,\C )$ as an additiv group and have for 
$Y(T)=\Hom (\C^* ,T)$ the natural perfect pairing $\langle-,-\rangle \colon X(T) \times Y(T) \to \Z, (\la, \sigma) \mapsto m$ with 
$ \la \circ \sigma (z) = z^m, z\in \C^*$. For the definition of the dual roots $\al_s^*\in Y(T), s\in S$ see \cite{CG}, chapter 7.1, p.361.     

\begin{satz} (\cite{CG}, thm 7.2.5, thm 8.1.16 - DL correspondence for affine Hecke algebras)\\
Let $G$ be a connected, simply connected semi-simple group over $\C$. 
\begin{itemize}
\item[(a)]
It holds $K_0^{G\times \C^*}(Z)\cong \mcH$ where $\mcH$ is the affine Hecke algebra associated to $(W,S)$, i.e.,
the $\Z[q,q^{-1}]$-algebra generated by $\{ e^{\la}T_w \mid w\in W, \la \in X(T), e^0=1\}$ with relations 
\begin{itemize}
\item[(i)] $(T_s+1)(T_s-q)=0, s\in S$, and $T_xT_y =T_{xy}$ for $x,y\in W$ with $\ell (xy)=\ell (x) \ell (y) $.
\item[(ii)] The $\Z[q,q^{-1}]$-subalgebra spanned by $e^{\la}$ is isomorphis to $(\Z [q^{\pm}])[X_1^{\pm}, \ldots, X_n^{\pm}]$, $n=rk (T)$.  
\item[(iii)] 
 For $\langle \la , \al_s^*\rangle =0$ it holds $T_se^{\la} =e^{\la}T_s$.\\
 For $\langle \la , \al_s^*\rangle =0$ it holds $T_se^{\la}T_s= qe^{\la} $. 
\end{itemize} 

\item[(b)]
The operation of $\mcH$ on a simple module factors over $H_*(Z^a)$, with $a=(s,t) \in G\times \C^*$ a semisimple element, in particular is $H_*((\pi^{-1}(x))^s)$ via the convolution construction a $\mcH$- module. The operation of the component group $C(a)=\Stab_{G\times \C^*}(a)/(\Stab_{G\times \C^*}(a))^o$ on $H_*(\pi^{-1}(x)^s)$ commutes with the $H_*(Z^a)$-operation and gives 
$H_*(\pi^{-1}(x)^s)=\bigoplus_{\chi \in \Simp(C(a))} K_{a,x,\chi}\otimes \chi $ for some $H_*(Z^a)$-modules $K_{a,x,\chi}$ (\textbf{ standard modules}).\\
If $t\in \C$ is not a root of unity, then there is a(n explicit) bijection between
\begin{itemize}
\item[(1)] $
\{ G-\text{conj. cl. of }(s,x,\chi)  \mid s\in G \text{ semisimple},\; sxs^{-1}=tx, \chi \in \Simp(C(s,t)), K_{(s,t),x,\chi}\neq 0 \}$
\\ and 
\item[(2)] Simple $\mcH$-modules where $q$ acts by multiplication with $t$. 
\end{itemize}
All simples are constructed from the standard modules, in general it is difficult to determine when the candidates are nonzero. 
For $t$ a root of unity there is an injection of the set (2) in (1). 
\end{itemize}
\end{satz}
%For the 
%Also include: Chriss-Ginzburg's result on affine Weyl groups and affine Hecke algebras (Deligne-Langlands-Lusztig correspondence ??) .
%\input{01-chapter-qgSpringerCorr}

\subsection*{Quiver-graded Springer Theory}\label{SpringerCor}

Let $Q$ be a finite quiver with set of vertices $Q_0$ and set of arrows $Q_1$. Let us fix a dimension vector $\dd \in \N_0^{ Q_0}$ and a sequence 
of dimension vectors $\ddd := (0=\dd^0, \ldots , \dd^{\nu}=:\dd), \dd^k_i\leq \dd_i^{k+1}$.  
Quiver-graded Springer Theory arises from the following initial data 
\[
\left[
\begin{aligned}
(*)&\; G=\Gl_{\dd}:=\prod_{i\in Q_0} \Gl_{d_i},\\ 
(*)&\; P= P(\ddd):= \prod_{i\in Q_0} P(d_i^{\bullet})\text{ where }P(d_i^{\bullet})\text{ is the parabolic in }\Gl_{d_i}\text{ fixing a (standard) flag} \\ &\; V_i^{\bullet}\text{ in }\C^{d_i}\text{ with dimensions given by }d_i^{\bullet},\\ 
(*)&\;  V=\Rd:= \prod_{(i\to j)\in Q_1} \Hom (\C^{d_i}, \C^{d_j})\text{ with the operation }(g_i) (M_{i\to j}) = (g_j M_{i\to j}g_i^{-1})\\
&\; \text{ is called representation space.}\\ 
(*)&\; F= F(\ddd ):=\{ (M_{i\to j}) \in \Rd \mid M_{i\to j}(V_i^k)\subset V_j^k , \;\; 0\leq k\leq \nu \}  
\end{aligned}
\right.
\]

Given $\dd$ and an (arbitrary) finite set $I:=\{ \ddd =(0=\dd^0, \ldots , \dd^{\nu})\mid \nu\in \N,  \dd^{\nu} =\dd\}$, we can describe the quiver-graded Springer 
correspondence explicitly via for $\ddd\in I$
\[  E_{\ddd} = 
\{ (M, U^\bullet) \in \Rd \times \Fl_{\ddd}  \mid \; i\xrightarrow{\forall \al} j \in Q_1\colon   M_{\al }(U_i^k)\subset U_j^k,\;  1\leq k \leq \nu\}
\]
\[ 
\xymatrix{
& E_{\ddd} \ar[dl]_{pr_1}\ar[dr]^{pr_2}  &\\
%\{ (M, U) \in \Rd \times \Fl_{\ddd}  \mid i\xrightarrow{\forall \al} j \in Q_1,  M_{\al }(U_i^k)\subset U_j^k,  1\leq k \leq \nu\}
\Rd && \Fl_{\ddd}}
\] 
where $\Fl_{\ddd} =\prod_{i\in Q_0} \Fl_{\ddd_i}$ and $\Fl_{\ddd_i}$ is the variety of flags of dimension $(0,\dd_i^1,\dd_i^2,\ldots , \dd_i^{\nu}=\dd_i)$ 
inside $\C^{\dd_i}$ and we set $E:= \bigsqcup_{\ddd \in I}E_{\ddd}$,  
\[ 
\begin{aligned}
Z_{\ddd, \ddd^{\prime}} &:= E_{\ddd}\times_{\Rd} E_{\ddd^\prime} \\
                        &\{ (M, U^\bullet,V^\bullet) \in \Rd \times \Fl_{\ddd}\times \Fl_{\ddd^\prime}\mid  i\xrightarrow{\forall \al} j \in Q_1 \colon   M_{\al }(U_i^k)\subset U_j^k, M_{\al }(V_i^k)\subset V_j^k \}
                        \end{aligned}
                        \]
\[
\xymatrix{ &Z_{\ddd, \ddd^{\prime}} \ar[dl]_{pr_1}\ar[dr]^{pr_{2,3}}  &\\
\Rd && \Fl_{\ddd}\times \Fl_{\ddd^\prime}}
\]                        
and the Steinberg variety is $Z:= \bigsqcup_{\ddd, \ddd^{\prime}\in I}  Z_{\ddd, \ddd^{\prime}}$.                        
This description goes back to Lusztig (cp. for example \cite{L}). It holds 
\[\dim E_{\ddd} =\dim \Fl_{\ddd} + \dim F( \ddd )= \sum_{i\in Q_0} \sum_{k=1}^{\nu -1}d_i^k(d_i^{k+1}-d_i^k) + \sum_{(i\to j)\in Q_1}\sum_{k=1}^{\nu} (d_i^k-d_i^{k-1})d_j^k, \] 
We define $\tits{\ddd}{\ddd}:= \dim \Gd-\dim E_{\ddd}$ and when $Q$ is without oriented cycles this is the tits form for the algebra 
$\C Q\otimes \C\mathbb{A}_{\nu +1}$ (cp. \cite{W}, Appendix)  
\[\tits{\ddd}{\ddd} =\sum_{k=0}^{\nu}\tits{\dd^k}{\dd^k}_{\C Q}-\sum_{k=0}^{\nu-1} \tits{\dd^k}{\dd^{k+1}}_{\C Q}.\] 
  
Let us take $(\ddd_i)_{i\in I}$ be the set of complete dimension filtrations of a given dimension vector $\dd$. The ($\Gd$-equivariant) Steinberg 
algebra is the quiver Hecke algebra (for $Q,\dd $). If the quiver $Q$ has no loops, the image of the injective map from lemma \ref{inj} has been 
calculated by Varagnolo and Vasserot in \cite{VV}. With generators and relations of the algebra they check that this is the same algebra as has 
been introduced by Khovanov and Lauda in \cite{KL} (and which was previously conjectured by Khovanov and Lauda to 
be the Steinberg algebra for quiver-graded Springer 
theory with complete dimension vectors). 
Independently, this has been proven by Rouquier in \cite{Ro}. 
\begin{satz} (quiver Hecke algebra, \cite{VV}, \cite{Ro}) Let $Q$ be a quiver without loops and $\dd \in \N_0^{Q_0}$ be a fixed dimension vector. 
The ($\Gd$-equivariant) 
quiver-graded Steinberg algebra for complete dimension filtrations $R_{\dd}^G:=H_*^G(Z)$ for $(Q, \dd)$ is as graded 
$\C$-algebra generated by 
\[ 1_i, i\in I, \quad z_i(k), i\in I, 1\leq k \leq d, \quad \si_i(s) , i\in I, s\in \{ (1,2), (2,3),\ldots , (d-1,d)\}=:\S, 
 \]
where $d:= \sum_{a\in Q_0}d_a$, $I:= I_{\dd}:= \{ (i_1,\ldots , i_d)  \mid \; i_k \in Q_0,  \sum_{k=1}^d i_k =\dd \}$ and we see $\S\subset S_d$ as permutations of $\{1,\ldots ,d\}$, we also define 
\[h_i((\ell, \ell+1)) = h_{i_{\ell +1}, i_{\ell}} =\# \{ \al \in Q_1\mid \al\colon i_{\ell +1} \to i_{\ell}\}\]
and let 
\[ \deg 1_i =0, \; \; \deg z_i (k) = 2, \; \; \deg \si_i((\ell, \ell +1)) =\begin{cases}
2h_i((\ell, \ell+1))-2 &, \text{ if } i_{\ell}=i_{\ell +1} \\
2h_i((\ell, \ell +1)) &, \text{ if }i_{\ell}\neq i_{\ell +1}  
\end{cases}  
\]
\\  
subject to  relations
\begin{itemize}

\item[(1)] (\emph{orthogonal idempotents})
\[
\begin{aligned}
1_i 1_j  &=\delta_{i,j} 1_i,\\
 1_i \si_i(s) 1_{is} &= \si_i(s) \\
 1_i z_i(k) 1_i &= z_i(k) 
\end{aligned}
\]

\item[(2)] (\emph{polynomial subalgebras}) 
\[z_i(k) z_i(k^\prime)= z_i(k^\prime) z_i(k)\]

\item[(3)]
For $s=(k,k+1)$, $i=(i_1,\ldots , i_d)$ we write $is:= (i_1,\ldots , i_{k+1}, i_k, \ldots , i_d)$ and set 
\[ \al_s := \al_{i,s}:= z_i(k)- z_i(k+1)\]
if it is clear from the context which $i$ is meant. We denote by $h_i(s):= \#\{ \al \in Q_1\mid \al\colon i_k \to i_{k+1}\}$.  
\[
\si_{i}(s) \si_{is}(s)= \begin{cases} 0 , &\text{ if } is=i\\
 (-1)^{h_{is}(s)} \al_s^{h_i(s)+h_{is}(s)} , & \text{ if } is\neq i.
 \end{cases}
 \]  
\item[(4)] (\emph{straightening rule})\\
For $s=(\ell , \ell +1)$ we set 
\[s(z_i(k))= \begin{cases} z_i(k+1), &\text{ if } k=\ell \\
z_i(k-1)  &\text{ if } k=\ell+1 \\
z_i(k)&\text{ else. }
\end{cases}
\]
\[ 
\si_i(s) z_{is}(k) - s(z_{is}(k))\si_i(s) = \begin{cases} -1_i, &\text{ if } is=i, s=(k,k+1) \\
1_i &\text{ if }is=i, s=(k-1,k) \\
0 &\text{, if } is\neq i. 
\end{cases}
\]
\item[(5)] (\emph{braid relation})\\
Let $s,t \in \S, st=ts$, then 
\[ \si_i(s)\si_{is}(t)  = \si_i (t)\si_{it }(s) . 
\]
 Let $i\in I$, $s=(k,k+1), t=(k+1, k+2)$. We set $s(\al_t):= (z_i(k)-z_i(k+2))=:t(\al_s)$ 
\[ 
\si_{i}(s)\si_{is}(t)\si_{ist} (s) -\si_{i}(t)\si_{it}(s)\si_{its} (t) = \begin{cases}
 P_{s,t} &\text{ if } ists =i , is\neq i, it \neq i \\
 0 , &\text{ else.}
\end{cases}
 \] 
where \[P_{s,t}:= \al_s^{h_{i}(s)} \frac{\al_t^{h_{is}(s)} - (-1)^{h_{is}(s)}\al_s^{h_{is}(s)}}{\al_s+\al_t} - \al_t^{h_{is}(s)} \frac{\al_s^{h_{i}(s)}-(-1)^{h_{i}(s)}\al_t^{h_{i}(s)}}{\al_s+\al_t} \] 
 is a polynomial in $z_i(k), z_i(k+1), z_i(k+2)$. 
\end{itemize}
We call this the \textbf{quiver Hecke algebra } for $Q, \dd$.   
\end{satz}

using the degeneration of the spectral sequence argument from lemma \ref{oddVan}  we get 
\begin{coro}\label{quivCox} Let $Q$ be a quiver without loops and $\dd\in \N_0^{Q_0}$. The not-equivariant Steinberg algebra $R_{\dd}:=H_{[*]}(Z)$ is 
the graded $\C$-algebra generated by 
\[ 1_i, i\in I, \quad z_i(k), i\in I, 1\leq k\leq d \quad \si_i(s) , i\in I, s\in \{ (1,2), (2,3),\ldots , (d-1,d)\} \]
with the same degrees and relations as $R_{\dd}^G$ and the additional relations
\[ P (z_i(1), \ldots , z_i(n)) =0 ,\quad i \in I, \;\; P \in \C[x_1,\ldots x_d]^{S_d}.\]
  
% \begin{itemize}

% \item[(1)] (\emph{orthogonal idempotents})
% \[
% \begin{aligned}
% 1_i 1_j &=\delta_{i,j} 1_i,\\
% 1_i \si_i(s) 1_{is} &= \si_i(s)
% \end{aligned}
% \]

% \item[(2)]
% \[ \si_{i}(s) \si_{is}(s)= \begin{cases}1, &\text{ if } is\neq i, h_i(s)=h_{is}(s)=0 \\
% 0 , & \text{ else.}
% \end{cases}
% \]

% \item[(3)] (\emph{braid relation})\\
% Let $s,t\in \S, st=ts$, then 
% \[ 
%  \si_i(s)\si_{is}(t)  = \si_i (t)\si_{it }(s) . 
% \]
%  Let $s=(k,k+1), t=(k+1,k+2)\in \S$, $i\in I$. 
% \[ 
% \si_{i}(s)\si_{is}(t)\si_{ist} (s) -\si_{i}(t)\si_{it}(s)\si_{its} (t) = \begin{cases}
%  -1_i, &\text{ if } ists =i , is\neq i \neq it,\\
%  & \quad \quad  h_i(s)+h_{is}(s)=1  \\
%  0 , &\text{ else.}
% \end{cases}
%  \]  
 
% \end{itemize}
%We name this algebra \textbf{quiver Coxeter algebra}\footnote{we are not aware that it has already a name}. 
\end{coro}

\paragraph{What about Springer fibre modules and the decomposition theorem?}

This is not investigated yet. We make some remarks on it. 
\begin{rem}
\begin{itemize}
\item[(1)] If $Q$ is a Dynkin quiver\footnote{i.e. the underlying graph is a Dynkin diagram of type $A_n,D_n , E_{6/7/8}$.}, 
the images of all quiver-graded Springer maps have finitely many orbits. %For a general quiver and dimension vector, it usually has 
%infinitely many orbits but of course there are always some (low dimensional) cases with finitely many orbits. 
For all quiver $Q$ and dimension vector $\dd \in \N_0^{Q_0}$ all $\Gl_{\dd}$-orbits in $\Rd$ are connected, 
i.e. $C(x)=\{e\}$ for all $x\in \Rd$. %If there are infinitely many orbits it is not known which stratification makes the quiver-graded Springer map locally trivial. 
\item[(2)] %There are very few cases in which the quiver-graded Springer maps are semi-small. 
In the case of finitely many orbits in the image of the Springer map, semi-smallness of the Springer map (associated to a dimension filtration $\ddd$ of a dimension vector $\dd$) 
is equivalent to for every $x\in \Rd$ it holds 
\[ 2 \dim \pi_{\ddd}^{-1}(x) \leq \dim \Ext^1_{\C Q} (x,x)= \codim_{\Rd} Gx. \]
It is very rarely fulfilled.  
\item[(3)] If $Q$ is a Dynkin quiver and $\dd\in \N_0^{Q_0}$ a complete set of the isomorphism classes of 
simple modules for the quiver Coxeter algebra $R_{\dd}$ is parametrized by the $G:=\Gl_{\dd}$-orbits in $\Rd$. For $x\in \Rd$ we have a simple
module of the form 
\[ L_{Gx} := \bigoplus_{\ddd } L_{Gx}^{(\ddd)}\]
where $\ddd$ runs over all complete dimension filtrations of $\dd$ and $L_{Gx}^{(\ddd)}$ is the multiplicity vector space occurring in the decomposition of 
$(\pi_{\ddd})_*\underline{\C}[e_{\ddd}]$. 
By the work of Reineke (see \cite{R}) there exists for every $x\in \Rd$ a complete 
dimension filtration $\ddd$ such that $Gx$ is dense in the image of $\pi_{\ddd}$. This implies by the considerations on page 11 that 
\[
L_{Gx, -*}^{(\ddd)}(-d_{Gx}) = H_{[*]}(\pi_{\ddd}^{-1}(x)) \quad (\neq 0),
\]  
as graded vector spaces, where $d_{Gx}= \dim Gx$. In fact, Reineke even shows that there exists a $\ddd$ for every $x$ such that the Springer map is a bijection over $Gx$, in which case $\dim L_{Gx}^{(\ddd)} =1$. 
\end{itemize}
\end{rem}

For $Q$ Dynkin, 
there are parametrizations of indecomposable graded projective modules in terms of Lyndon words, see \cite{HMM},
which are not yet understood in the context of the decomposition theorem. 

\subsubsection*{Monoidal categorifications of the negative half of the quantum group}

Again let $Q$ be a finite quiver without loops. 
First Lusztig found the monoidal categorification of the negative half of the quantum group via perverse sheaves, then Khovanov and Lauda did the same with (f.g. graded) projective modules over quiver Hecke algebras. In the following theorem's we are explaining the following diagrams of isomorphisms of twisted Hopf 
algebras over $\Q (q)$.  
\[
\xymatrix{
&K_0(proj^{\Z} \bigoplus H_*^G(Z))\otimes \Q(q) \ar[dl]\ar[dr] & \\
\mcU^-:=U_q^-(Q) & & K_0(\mcP )\otimes \Q (q) \ar[ll]
}
\]
In all three algebras there exists a notion of canonical basis which is mapped to each other under the isomorphisms. 
Also, there is a triangle diagram with isomorphisms defined over $\Z[q,q^{-1}]$ which gives the above situation after applying 
$-\otimes_{\Z[q,q^{-1}]} \Q (q)$. 

\paragraph{The negative half of the quantum group.}
The negative half $\mcU^-:=U_q^-(Q)$ of the quantized enveloping algebra (defined by Drinfeld and Jimbo) associated to the quiver $Q$ is defined via: 
Let $a_{i,j}:= \# \{ \al\in Q_1\mid \al \colon i\to j, \text{ or }\al\colon j\to i\}, i\neq j \in Q_0$. 
It is the $\Q(q)$-algebra generated by $F_i, i\in Q_0$ with respect to the (quantum Serre relations) 
\[ 
\sum_{p=0}^{N+1} [p, N+1-p] F_j^pF_i F_j^{N+1-p} =0 ,\quad  N=a_{ij}, i\neq j
\]
where \[ [n]_!:= \prod_{k=1}^n \frac{q^k-q^{-k} }{q-q^{-1}}, \quad [n,m]:= \frac{[n+m]_!}{[n]_![m]_!}. \]
Lusztig calls this ${}^{\prime}\mathbf{f}$

A Hopf algebra is a bialgebra (i.e. an algebra which also has the structure of a coalgebra such that the comultiplication and counit are algebra homomorphisms) 
which also has an antipode, i.e. an anti-automorphism which is uniquely determined by the 
bialgebra through commuting diagrams. 
A twisted Hopf algebra differs from the Hopf algebra by: The comultiplication and the antipode are only homomorphisms if you \emph{twist} the algebra 
structure by a bilinear form (see the example below). For more details on the definition see \cite{LZ}.  
The twisted $\Q (q)$-Hopf algebra structure  is given by the following, it is by definition a $\Q (q)$-algebra which is $\N^{Q_0}$-graded and it has  
\begin{itemize}
\item[(1)] (comultiplication) \\
If we give $\mcU^{-}\otimes_{\Q(q)}\mcU^{-}$ the algebra structure 
\[ (x_1\otimes x_2)(x_1^{\prime}\otimes x_2^{\prime} := q^{\left| x_2 \right| \cdot \left|x_1^{\prime}\right|} x_1 x_1^{\prime}\otimes x_2x_2^{\prime}\] 
where for $x\in \mcU^{-}$ we write $\left| x\right|\in \N_0^{Q_0} $ for its degree and the symmetric bilinear form 
\[ \cdot \colon \Z_0^{Q_0} \times \Z^{Q_0} \to \Z, \;\;i\cdot i := 2,\;  i\cdot j := -a_{i,j} \text{ for }i\neq j \]
Then the comultiplication is the $\Q (q)$ algebra homomorphism 
\[
\mcU^{-} \to \mcU^{-} \otimes \mcU^{-}, \; \; F_i \mapsto F_i\otimes 1 + 1\otimes F_i
\]
\item[(2)] (counit) $ \; \; \epsilon \colon \mcU^{-} \to \Q (q) , \;\; F_i \mapsto 0 $
\item[(3)] (antipode)\\
Let $\mcU^{-}_{tw}$ be the algebra with the multiplication $x * y:= q^{\left|y \right| \cdot \left|x \right|}xy$ \\
The antipode is the algebra anti-homomorphism 
\[\mcU^{-}\to \mcU^{-}_{tw}, \; \; F_i\mapsto -F_i \]
\end{itemize}

\paragraph{Lusztig's category of perverse sheaves.}

Lusztig writes complete dimension filtrations as words in the vertices $i=(i_1,\ldots , i_d)$, $i_t\in Q_0$, set $\dd:= \sum_{t=1}^di_t$ 
and defines 
\[L_i:= (\pi_i)_* \underline{\C}[e_i]\] 
where $\pi_i \colon E_i:= \Gl_{\dd}\times^{P_i}F_i\to  R_Q(\dd )$ is the quiver-graded Springer map and $e_i=\dim_{\C} E_i$.  
Let us call $\mcP_{Q_0}$  the additive category generated by shifts of the $L_i, \; i=(i_1,\ldots , i_d), i_t\in Q_0$. 
The homomorphisms $\Hom (L_i, L_j[n])$ in this category is zero unless $\dd=\sum i_t =\sum j_k$ and then they are given by 
$1_j* H_{[n]}^{\Gl_{\dd}} (Z) *1_i$. It can be endowed with the structure of a monoidal category via 
\[ 
L_i*L_j:= L_{ij}
\]
where $ij$ is the concatenation of the sequence $i$ and then $j$.

\begin{lemma} (Lusztig, \cite{L}, Prop. 7.3)
Let $\mcP$ be the idempotent completion of $\mcP_{Q_0}$ (i.e. we take the smallest additive category 
generated by direct summands of the $L_i$ in $D_{\Gl_{\dd}}^b(R_Q(\dd))$ and their shifts). It carries a monoidal structure and the inclusion induces 
\[
K_0(\mcP_{Q_0}) = K_0 (\mcP)
\]
where the Grothendieck group has the ring structure from the monoidal categories and a $\Z [q,q^{-1}]$-module structure via the shift, i.e. 
$q\cdot[M]:= [M[1]]$, $M$ an object in $\mcP$. 
\end{lemma}
\begin{rem}
We call the monoidal category $\mcP$ \textbf{Lusztig's category of perverse sheaves}. Even though these are not perverse sheaves since we allow shifts of them and 
Lusztig originally defined them inside $\bigsqcup D^b(R_Q(\dd ))$ which of course gives a different category (for example in this category 
$\Hom (L_i, L_j[n] ) = 1_j *H_{[n]}(Z)*1_i$). 
Nevertheless the two categories have the same Grothendieck group. 
In the view of the context here we think it is more 
apropriate to define it in the equivariant derived categories.  
\end{rem}

\begin{rem}
The previous lemma is no longer true if you allow your quiver to have loops. \\
For example if $Q$ is the quiver with one loop. Then, let $Z_n$ be the Steinberg algebra associated to $(G=\Gl_n,B_n, \mathfrak{g}\mathfrak{l}_n, \mathfrak{n}_n)$ 
with $B_n\subset \Gl_n$ the upper triangular matrices, $\mathfrak{n}_n$ the Lie algebra of the unipotent radical of $B_n$. We claim  
\[ \begin{aligned} K_0(\mcP) & = \bigoplus_{n\in \N_0}K_0( \text{f.d. proj. graded }H_{[*]}(Z_n)-\text{modules})\\
                            & = \bigoplus_{n\in \N_0} K_0(\text{f.d. simple graded }H_{[0]}(Z_n)-\text{modules} )\\
                            & = \bigoplus_{n\in \N_0} K_0(\text{f.d. graded }\C S_n-\text{modules} ) =  
\left(\bigoplus_{n\in \N_0} K_0(\text{f.d.}\C S_n-\text{modules} )\right) \otimes _{\Z} \Z[q,q^{-1}] \\
& = (\text{Symmetric functions}) \otimes_{\Z} \Z [q,q^{-1}] 
\end{aligned}
\]  
The first isomorphism is implied by the Corollary \ref{projectives}. The second equality is implied by semi-smallness of the classical Springer maps.
For the third result see the section on classical Springer theory. The last equality is well-known, it maps the simple module $S_{\la}$ (=Specht module) corresponding to a partition $\la$ to the Schur function corresponding to $\la$. \\
But the category $\mcP_{Q_0}$ corresponds to the submonoidal category given by finite direct sums of shifts of 
finite-dimensional free modules. This is a monoidal category generated by direct sums of shifts of  
one object $E=S_1$ and an arrow $s\colon E^2:=E\otimes E \to E^2$  of degree $0$ with the 
relation $(sE)\circ (Es)\circ (sE)=(Es)\circ (sE)\circ (Es)$ (see also \cite{Ro}). 
In this case $K_0(\mcP_{Q_0})= \Z[q,q^{-1}, T], \; [E] \mapsto T$ which is much smaller than $K_0(\mcP)$.    
\end{rem}

Now, let again be $Q$ without loops. 
$K_0(\mcP)$ has the structure of a twisted $\Z[ q, q^{-1}]$-Hopf algebra. The algebra structure is given by the 
monoidal structure on $\mcP$ which is defined by induction functors. 
A restriction functor for the category $\mcP$ defines the  
structure of a coalgebra. For the geometric construction of these functors see \cite{L}. 
% \begin{itemize}
% \item[(1)] (comultiplication)
% \item[(2)] (counit)
% \item[(3)] (antipode)
% \end{itemize}

\begin{satz} (Lusztig, \cite{L}, thm 10.17)
The map 
\[ 
\begin{aligned}
\la_Q\colon \mcU^- &\to K_0(\mcP)\otimes_{\Z[q,q^{-1}]} \Q (q)\\
F_i & \mapsto [L_i] \otimes 1, \quad i\in Q_0
\end{aligned}
\]
here we see $i\in Q_0$ as a sequence in the vertices of length $1$. This defines an isomorphism of twisted $\Q (q)$-Hopf algebras. 
\end{satz}

\begin{defi}
We call $\mathbf{B}:=\{ [L_i]\otimes 1 \mid i=(i_1,\ldots , i_d), i_t\in Q_0\}$  
\textbf{canonical basis} for $K_0(\mcP )\otimes \Q (q)$. \\
We also call $\la_Q^{-1}(\mathbf{B})$  \textbf{canonical basis} in $\mcU^{-}$. \\
Also the image in $K_0(proj^{\Z} \bigoplus R_{\dd}^{\Gl_{\dd}} )\otimes \Q (q)$ is called  \textbf{canonical basis}. 
\end{defi}

There are two intrinsic alternative definitions of the canonical basis for $\mcU^{-}$ given by again 
Lusztig in \cite{L3} for the finite type case and in general by Kashiwara's crystal basis, see \cite{Kash}.

\paragraph{Generators and relations for $\mcP_{Q_0}$.}
This is due to Rouquier (cp. \cite{Ro}), it is the observation that the generators and relations of the quiver Hecke algebra rather easily give 
generators and relations for the monoidal category $\mcP_{Q_0}$. In the category, we use the convention instead of $E\to E^\prime (n)$ 
we write $E\to E^\prime$ is a morphism of degree $n$. A composition $g\circ f$ of a morphism $f\colon E\to E^\prime$ of degree $n$ and 
$g\colon E^\prime \to E^{\prime \prime}$ of degree $m$ is the homomorphism $E\to E^{\prime \prime}$ of degree $n+m$ given by 
$E\xrightarrow{f} E^\prime(n) \xrightarrow{g(n)}E^{\prime \prime}(n+m)$.  \\

Let $Q$ be a quiver without loops.
Let $\mcB$ be the monoidal category generated by finite direct sums of shifts of objects $E_a =: E_a(0), a\in Q_0$ and arrows 
\[ z_a\colon E_a\to E_a,\;\; \si_{a,b} \colon E_aE_b \to E_bE_a, \;\;  a,b \in Q_0\]
of degrees 
\[\deg z_a =2,\;\;  \deg \si_{a,b}= \begin{cases}-2 &, \text{ if }a=b\\
2h_{b,a} &, \text{ if }a\neq b 
\end{cases}
\]
where as before $h_{a,b}:= \# \{ \al \in Q_1 \mid \al\colon a\to b\}, \;\; a,b \in Q_0.$ and assume relations 
\begin{itemize}
\item[(1)] ($s^2=1$)
\[ 
\begin{aligned}
\si_{ab}\circ \si_{ba} &= \begin{cases} 
(-1)^{h_{b,a}}(E_bz_a-z_bE_a)^{h_{a,b}+h_{b,a}}&,\;\; \text{if } a\neq b \\
0 &,\;\; \text{if } a=b \\
\end{cases}\\
\end{aligned}
\]
\item[(2)] (straightening rule)
\[
\begin{aligned}
\si_{ab}\circ z_aE_b-E_bz_a\circ \si_{ab} & =
\begin{cases}
0 &\quad\! ,\;\; \text{if } a\neq b, \\
E_aE_a&\quad\! ,\;\; \text{if } a= b, 
\end{cases} \\ 
\si_{ab}\circ E_az_b - z_bE_a\circ \si_{ab} &= 
\begin{cases}
0 &,\;\; \text{if } a\neq b, \\
-E_aE_a&,\;\; \text{if } a= b, 
\end{cases}\\
\end{aligned}
\]

\item[(3)] (braid relations) for $a,b,c\in Q_0$ we have the following inclusion of $\C$-algebras. Let $\C[\al_s,\al_t]$ be the set of polynomials in $\al_s, \al_t$. \[
\begin{aligned}
J_{a,b,c}\colon \C[\al_s, \al_t] & \to \End_{\mcB} (E_aE_bE_c) \\
\al_s & \mapsto E_az_bE_c -z_aE_bE_c \\
\al_t &\mapsto E_aE_bz_c - E_az_bE_c,
\end{aligned}
\]
we set $ t(\al_s^h):=(\al_s+ \al_t)^h =: s(\al_t^h)\in \C[\al_s, \al_t], h\in \N_0$. Then, the relation is  
\[
\begin{aligned}
\si_{ab}&E_c\circ E_a\si_{cb} \circ \si_{ca}E_b - E_b \si_{ca} \circ \si_{cb}E_a\circ E_c\si_{ab} \\
&= \begin{cases}
J_{bab}(\al_s^{h_{a,b}}\frac{\al_t^{h_{b,a}}- (-1)^{h_{b,a}}\al_s^{h_{b,a}}}{\al_s+\al_t} -\al_t^{h_{b,a}} \frac{\al_s^{h_{a,b}}-(-1)^{h_{a,b}}\al_t^{h_{a,b}}}{\al_s+\al_t} )  
&, \;\; \text{if } a=c, a\neq b,  \\
0 &, \;\; \text{else. }
\end{cases}
\end{aligned}
\]

\end{itemize}

for $i=(i_1,\ldots , i_n), i_t \in Q_0$ we set $E_i:= E_{i_1}E_{i_2}\cdots E_{i_n}$.
Let $I_{\dd}:= \{ i=(i_1,\ldots , i_n)\mid \sum_t i_t = \dd\}$.
Then, by construction there is an isomorphism of algebras
\[
\begin{aligned}
R_{\dd}^{\Gl_{\dd}} & \to \bigoplus_{i,j\in I_{\dd }} \Hom_{\mcB}(E_i, E_j)\\
1_i &\mapsto \id_{E_i} \\
z_i(t) &\mapsto E_{i_1}E_{i_{2}}\cdots E_{i_{t-1}} z_{i_t} E_{i_{t+1}} \cdots E_{i_n}\\
\si_i(s) &\mapsto E_{i_1}\cdots E_{i_{\ell-1}} \si_{i_{\ell+1}, i_{\ell}} E_{i_{\ell+2}}\cdots E_{i_n} , \quad \text{, if } s=(\ell,\ell+1)\in S_n
\end{aligned}
\]

 \begin{satz} (\cite{Ro}) There is an equivalence of monoidal categories 
\[
\begin{aligned}
 \mcP_{Q_0} & \to \mcB \\
 L_i &\mapsto E_i
 \end{aligned}
 \]
 which is on morphisms the isomorphism of algebras from above. 
 \end{satz}
%
%\begin{satz} (Rouquier, \cite{Ro}) Let $\mcC$ be the monoidal category with generators and relations ... 
%Then $\mcC \cong \mcP_{Q_0}$ as monoidal categories. 
%\end{satz}
Since we have not more knowledge on the decomposition theorem for quiver-graded Springer maps, we can not expect to find a similar 
easy description for the category $\mcP$. 

\paragraph{Khovanov and Lauda's catgorification of the negative half of the quantum group.}

Many years later Khovanov and Lauda have a different approach to the same \emph{monoidal categorification} as Lusztig. 
Instead the category $\mcP$ they consider the category of projective graded (f.g.) modules over quiver Hecke algebras $R_{\dd}:= R_{\dd}^{\Gl_{\dd}}$, $\dd \in \N_0^{Q_0}$ 
\[ 
proj^{\Z} \bigoplus_{\dd\in \N_0^{Q_0} } R_{\dd}
\]
It is easy to see that we have natural injective maps $\mu \colon R_{\dd}  \otimes R_{\ee }  \to R_{\dd + \ee} $ compatible with the 
algebra multiplication. We write $1_{\dd ,\ee}:= \mu (1\otimes 1)$. From this there are (well-defined see \cite{KL},section 2.6) 
induction and restriction functors 
\[
\begin{aligned}
\Ind_{\dd , \ee}^{\dd +\ee} \colon proj^{\Z}(R_{\dd}  \otimes  R_{\ee })  & \to proj^{\Z}(R_{\dd + \ee}), \quad X\mapsto R_{\dd +\ee} 1_{\dd , \ee} \otimes_{R_{\dd}\otimes R_{\ee}} X\\
\res_{\dd , \ee}^{\dd +\ee}\colon  proj^{\Z}(R_{\dd + \ee})               &\to  proj^{\Z}(R_{\dd}  \otimes  R_{\ee }), \quad Y\mapsto 1_{\dd ,\ee}Y 
\end{aligned}
\]
The induction functor gives $proj^{\Z} \bigoplus_{\dd\in \N_0^{Q_0} } R_{\dd}$ the structure of a monoidal category via $X\circ X^\prime := 
\Ind_{\dd ,\ee}^{\dd +\ee} X\boxtimes X^\prime$ where $X\boxtimes X^\prime$ is the natural graded $R_{\dd}\otimes R_{\ee}$-module structure.   \\ 
The twisted $\Z[q,q^{-1}]$-Hopf algebra structure on $K_0 (proj^{\Z} \bigoplus_{\dd\in \N_0^{Q_0} } R_{\dd}^{\Gl_{\dd}})$ is given by:\\
Obviously, it is a $\Z[ q, q^{-1}]$-algebra with $q$ operating as the shift $(1)$ on the graded modules, i.e. $q\cdot [M]:= [M(1)]$. 
The comultiplication is given by $[\res][P]:= \sum_{\dd, \ee \colon \dd +\ee =\ff }[\Res_{\dd , \ee}^{\ff} (P)]$. 

\begin{satz} (Khovanov, Lauda, \cite{KL})
The map 
\[
\begin{aligned}
\kappa_Q\colon \mcU^{-} &\to K_0(proj^{\Z} \bigoplus_{\dd\in \N_0^{Q_0} } R_{\dd}^{\Gl_{\dd}}) \otimes \Q (q) \\
 F_i &\mapsto [R_i^{\Gl_1}] \otimes 1, \quad i\in Q_0
 \end{aligned}
\]
where we consider $i\in Q_0$ as an element in $\N_0^{Q_0}$, is an isomorphism of twisted $\Q (q)$-Hopf algebras.  
\end{satz}

We want to point out: Khovanov and Lauda invented the quiver Hecke algebra, which later turned out (by Varagnolo and Vasserot's result) to be 
the same as the Steinberg algebra of quiver-graded Springer theory. Their explicit description (generators and relations for the algebra) 
and diagram calculus (which we leave out in this survey) 
are a major step forward from Lusztig's description. Their work sparked a big interest in this subject.   

\begin{rem} Let $Q$ be a Dynkin quiver. Then, the objects of the category $\mcP$ are direct sums of shifts of 
$IC_{\mcO}$ where $\mcO \subset R_{Q}(\dd )$ is a $\Gl_{\dd}$-orbit (we do not write a local system if the trivial local system is meant).
These are in bijection with isomorphism classes of $\C Q$-modules. The monoidal structure on $\mcP$ is constructed such 
that $ K_0(\mcP)\otimes \Q (q)$ is the twisted Ringel-Hall algebra (over $\Q(q)$). \\
The isomorphism between the twisted Ringel-Hall algebra and the negative half of the quantum group associated to the underlying graph of the quiver 
is a theorem of Ringel, see for example \cite{Rin}.  
\end{rem}

\begin{landscape}
\begin{table}
\caption{List of known Steinberg algebras.} 
%For any algebraic group, let $H_*^H(-,\C)$ be $H$-equivariant Borel-Moore homology, recall that for a Steinberg variety $Z$, there is a convolution product on 
%$H_*^H(Z,\C)$ and on $K_0^H(Z)$, see \cite{CG}. We write $H_{top}(Z,\C)$ for the highest degree Borel Moore homology group, it is a subalgebra of $H_*(Z,\C)$, see loc.cit. ... Recall that it holds $H_*(Z,\C)=K_0(Z)$ (also as rings?)
\begin{center}
\begin{tabular}{l||c|c|c|c}%{ p{3cm}|| p{2cm} | p{3cm} | p{3cm} | p{2cm}} 
   &$H_{top}(Z,\C)$ & $H_*(Z,\C)$ & $H_*^G(Z,\C)$& $K_0^{G\times \C^*} (Z)\otimes_{\Z}\C$   \\
\hline 
                          
$(G,B,\mathfrak{g}, \mathfrak{u})$ & $\C W$  & $\C[\mathfrak{t}]/I_W\# \C[W]$  &  $\C[\mathfrak{t}]\# \C[W]$ & affine    \\
classical ST & && degenerate affine & Hecke algebra \\
&&& Hecke algebra & \\
\hline

$(G,B, \{0\}, \{0\})$  & $\C$ & $\End_{\C-lin}(H^*(G/B))$ &  $\End_{H_G^*(pt)}(H^*_G(G/B))$  & \\
 nil ST   &&  &(affine) nil Hecke & ? \\
 i.e. $Z=G/B\times G/B$   &&& algebra &\\
\hline

quiver-graded ST & ? & $R_{\dd}$ & quiver Hecke algebra & ?\\
(complete dim filtrations) &&  & ($=$ KLR-algebra)&
\end{tabular}
\end{center}
\end{table}

Further known examples are: 
\begin{itemize}
\item[(1)] There is an exotic Springer theory (by Kato \cite{Ka1}, \cite{Ka2}, Achar and Henderson \cite{He}). The Steinberg algebra 
$K_0^{G\times (\C^*)^3} (Z)\otimes_{\Z}\C$ is isomorphic to the Hecke algebra with unequal parameters of type $C_n^{(1)}$. Also Kato gave an exotic Deligne-Langlands correspondence. 
\item[(2)] Quiver-graded Springer theory for the oriented cycle quiver (allowing only nilpotent representations) gives that $H_*^G(Z)$ is isomorphic to the quiver Schur algebra (compare the work of Stroppel and Webster, \cite{SW}.) 
\end{itemize}

\end{landscape}

\begin{landscape}
\begin{figure}
 \begin{center}
  \begin{tikzpicture}[node distance=3cm, auto, >=stealth]
   % nodes
   
   \node[decision,text width=4cm]
                (c)          {\begin{tabular}{r|l}
                  &\\ 
                clas-  $\; $ & $\; $ quiver- \\
              sical  $\; $ & $\;$ graded \\
              Springer     $\; $ &$\;$  Springer \\
              Theory  $\; $ & $\;$ Theory \\ 
              & 
              \end{tabular}};

   \node[block] (d)  [right of=c, node distance=8.8cm]           {$U_q^-$ (ass. to quiver), canonical basis};
   \node[block] (h)  [left of=d, node distance=3.5cm]     {perverse sheaves, \cite{L}};
   \node[block] (f)  [above of=c, node distance=5.5cm]             {(quantum) Schur Weyl theory ?};
   \node[block] (e)  [above right of=c, node distance=6cm]           {quiver Schur alg., \cite{SW}};
   \node[block] (g)  [right of=f, node distance=8.8cm]           {quiver Hecke alg., \cite{VV}, \cite{Ro} };
  
   \node[block] (j)  [below of=d]                              {categories of flags of submodules, \cite{Ri4}};
   \node[block] (i)  [left of=j, node distance=4cm]     {quiver flag varieties \cite{W}};
   \node[block] (k)  [below of=j, node distance=3cm]           {Cluster algebras, \cite{CC}};
   \node[block] (l)  [below of=i, node distance=3cm]           {Res of Sing for orbit closures, \cite{R}};

   \node[block] (o)  [left of=c, node distance=8.8cm]            {parame- trize simple modules};
   \node[block] (a)  [below of=o]                                {Res of ADE-surf. sing., Mac Kay corresp.};
   \node[block] (b)  [right of=a, node distance=3.5cm]          {Symplectic Res of Sing};  
   \node[block] (n)  [below of=a]                               {braid gr. operation on Coh(E)};
   \node[block] (m)  [right of=n, node distance=4cm]            {NCR associated to Springer map};
   \node[block] (p)  [left of=f, node distance=8.8cm]            {various Hecke alg. (and $\C W$)}; 
   \node[block] (q)  [above right of=o]                               {exotic ST, \cite{Ka1}};

   % edges
   \draw[->] (c) -- node[above] {slice down}(a);
   \draw[->] (c) -- (b);
   \draw[->] (c) -- (e);
   \draw[->] (c) -- (h);
   \draw[->] (c) -- (i);
   \draw[->] (c) -- (l);
   \draw[->] (c) -- (m);
   \draw[->] (c) -- node[above] {Steinberg alg.}(p);
   \draw[->] (c) -- (q);
   \draw[->] (c) -- node[above] {Springer fibres}(o);
   \draw[->] (c) -- node[below] {Steinberg alg.} (g);
  % \draw[->] (c.north) to [out=170,in=45] node[above] {no} (b.north);
   %\draw[->] (d.south) to [out=210,in=20] (e.north);
   \draw[->] (f) -- (p);
   \draw[->] (f) -- (g);
   \draw[->] (h) -- node[above] {$K_0$}(d);
   \draw[->] (m) -- node[above] {\cite{Be}}(n);
   \draw[->] (i) -- node[above] {Euler char}(k);
   \draw[->] (i) -- (j);
   \draw[->] (q) -- (p);
   \draw[->] (g) -- node[right] {$K_0(proj^{\Z}-)$}(d);
   \draw[->] (p) -- node[left] {$\begin{matrix} \text{Springer- /}\\
                                                 \text{DL corresp.}\end{matrix}$}(o);
  \end{tikzpicture}
  \caption{Springer Theory and related fields}
  \label{flowchart}
 \end{center}
\end{figure}
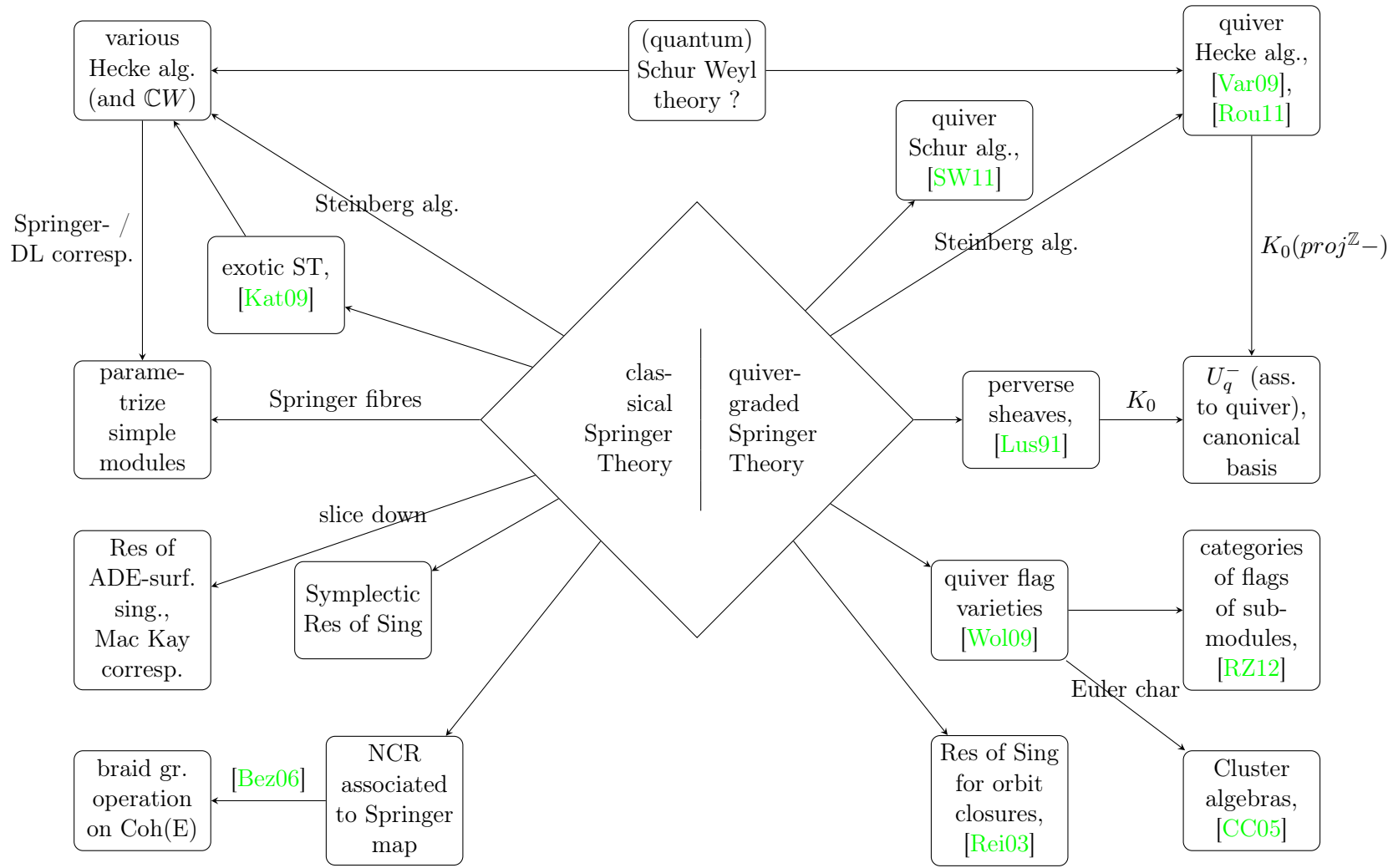

\end{landscape}

\subsection{Literature review}

Collapsings of homogeneous vector bundles are quite ubiquitous (for example see \cite{Ke}). %Later, Knutson and Shimonzono (\cite{KS} dropped the extra assumption and called the situation a \emph{Kempf collapsing}). 
\begin{itemize}
\item[(1)] Classical Springer Theory (cp remark \ref{others}):\\
Classical Springer Theory is usually defined for semi-simple algebraic groups and goes back to 
first Springer \cite{Spr1},\cite{Spr2}, then Kazhdan-Lusztig \cite{KL2}, Slodowy \cite{Sl2}, Lusztig \cite{L2}, Rossmann \cite{Ros}) 
and the defined convolution operations differ between each other by at most multiplication with a sign character (see \cite{Ho}). \\
Also relevant is the earlier work on the topology of Springer fibres of Spaltenstein (see \cite{Sp}, \cite{Sp2}) and Vargas (see \cite{V}) and the Springer map already occurrs in Steinberg's work (for example \cite{St2}). 
A book on classical Springer Theory is written by Borho, Brylinski and Mac Pherson \cite{BBM}.
A comprehensive treatment can be found in chapter 3 of \cite{CG} and a short one using perverse sheaves in 
\cite{Ar3} if you speak French. I apologize to the many other authors who I do not mention. 
%\item[(3)] Exotic Springer Theory:\\
%Kato \cite{Ka1}, \cite{Ka2}, (Achar and ) Henderson \cite{He1}
\item[(2)] Quiver-graded Springer Theory: \\
First considered by Lusztig, see \cite{L}. Later, Reineke started to look at it as an analogue of the 
classical Springer theory, see \cite{R}, also see \cite{W}. \\
The quiver Hecke algebras as Steinberg algebras first occurred in the work of Varagnolo and Vasserot, cp. \cite{VV}, and independently also 
in Rouquier's article \cite{Ro}. 
\end{itemize}

\paragraph{Open problems/ wild speculations: }
\begin{itemize}
\item[(O1)] Are Springer fibre modules always semi-simple modules over the Steinberg algebra?
\item[(O2)] Which Steinberg algebras are affine cellular algebras?\\
           Which have finite global dimension?\\
           Partial answers:  
           Brundan, Kleshchev and McNamara showed that KLR-algebras for Dynkin quivers are affine cellular (see \cite{BKM}).  \\
           Certain Steinberg algbras (including KLR-algebras for Dynkin quivers) have been shown to have finite global dimension (see \cite{Ka3}). 
           In \cite{BKM}, the authors write that they expect that KLR-algebras have finite global dimension if and only if the quiver is Dynkin. 
           
\item[(O3)] Are there Kazhdan-Lusztig polynomials and even a theory of canonical basis for Steinberg algebras?\\
           Do we have Standard modules for Steinberg algebras?\\
           Partial answers:
           Standard modules have been defined in \cite{Ka3} under some assumptions (finitely many orbits in the image of the Springer map,...). \\
           The original definition of Kazhdan-Lusztig polynomials has been inspired by studying a base change between two bases in the Steinberg variety                      associated to classical Springer theory. 
           
\item[(O4)] Can we describe noncommutative resolutions of singularities corresponding to Springer maps?\\ 
           Can we adapt the notion of a noncommutative resolution of singularities using constructible instead of coherent sheaves?\\
           Partial answers exists for the coherent sheaf theory: Bezrukavnikov studied it for classical Springer theory (see \cite{Be}) and 
           for quiver-graded Springer theory with $Q=A_2$ Buchweitz, Leuschke and van den Bergh studied noncommutative resolutions 
           (see \cite{BLvdB}, \cite{BLvdB2}). 
        
\item[(O5)] Does there exist a Schur-Weyl theory relating classical and quiver-graded Springer theory (for example via Morita equivalences of the associated 
           Steinberg algebras)? \\         
           Partial answers only for \emph{type A}-situations (so maybe it only exists in this case): due to Brundan, Kleshchev \cite{BKl}, see also for 
           example \cite{Web}.  
\end{itemize}

\addcontentsline{toc}{section}{\textbf{References} \hfill}

\bibliographystyle{alphadin}
\bibliography{SurveyOnSpringer}
\paragraph{Acknowledgements:} I would like to thank my supervisor Andrew Hubery for his constant support during my phd. 
Also, Bill Crawley Boevey for pointing out a mistake and Syu Kato for explaining me that the graded 
modules over $H_{<*>}^G (Z)$ and over $H_{[*]}^G(Z)$ are equivalent. 
I would also like to acknowledge financial support from the CRC 701 during my guest stays in Bielefeld.
\clearpage
\phantomsection
%\addcontentsline{toc}{section}{\textbf{Index} \hfill}
%\printindex

\end{document}